# BOOTSTRAP PERCOLATION IN THREE DIMENSIONS

By József Balogh,[1] Béla Bollobás[2] and Robert Morris[3]

*University of Illinois, University of Memphis and
University of Cambridge*

By bootstrap percolation we mean the following deterministic process on a graph $G$. Given a set $A$ of vertices "infected" at time 0, new vertices are subsequently infected, at each time step, if they have at least $r \in \mathbb{N}$ previously infected neighbors. When the set $A$ is chosen at random, the main aim is to determine the critical probability $p_c(G, r)$ at which percolation (infection of the entire graph) becomes likely to occur. This bootstrap process has been extensively studied on the $d$-dimensional grid $[n]^d$: with $2 \le r \le d$ fixed, it was proved by Cerf and Cirillo (for $d = r = 3$), and by Cerf and Manzo (in general), that

$$p_c([n]^d, r) = \Theta\left(\frac{1}{\log_{(r-1)} n}\right)^{d-r+1},$$

where $\log_{(r)}$ is an $r$-times iterated logarithm. However, the exact threshold function is only known in the case $d = r = 2$, where it was shown by Holroyd to be $(1 + o(1))\frac{\pi^2}{18 \log n}$. In this paper we shall determine the exact threshold in the crucial case $d = r = 3$, and lay the groundwork for solving the problem for all fixed $d$ and $r$.

**1. Introduction.** In this paper we shall study three-neighbor bootstrap percolation on $[n]^3$. Let $G$ be a (finite) graph, let $r \in \mathbb{N}$, and let $A \subset V(G)$ be a set of initially "infected" vertices. In $r$-neighbor bootstrap percolation on $G$, with initial set $A$, new vertices of $G$ are infected if they have at least $r$ infected neighbors, and infected vertices remain infected forever. Formally, set $A_0 = A$, and

$$A_{t+1} := A_t \cup \{v \in V(G) : |\Gamma(v) \cap A_t| \ge r\}$$

Received June 2008.
[1]Supported in part by OTKA Grant T049398, NSF Grant DMS-06-00303 and UIUC Campus Research Board 07048.
[2]Supported in part by ITR Grant CCR-0225610 and ARO Grant W911NF-06-1-0076.
[3]Supported in part by MCT Grant PCI EV-8C.
*AMS 2000 subject classifications.* 60K35, 60C05.
*Key words and phrases.* Bootstrap percolation, sharp threshold.

for each integer $t \geq 0$. The *closure* of $A \subset V(G)$ is the set $[A] = \bigcup_t A_t$ of eventually infected vertices. We say that the set $A$ *percolates* if eventually the entire vertex set is infected, that is, if $[A] = V(G)$.

Bootstrap percolation is an example of a cellular automaton, studied, for example, by von Neumann [28]. However, the particular model we are studying was introduced in 1979 by Chalupa, Leith and Reich [14], and was subsequently rediscovered by several authors who were motivated by its connections to interacting particle systems and other physical applications (see, e.g., the survey [1]). The first mathematical papers in the area were by van Enter [26] and Schonmann [24, 25], who studied the process on the infinite lattice $\mathbb{Z}^d$, and by Aizenman and Lebowitz [2], who studied it on the finite grid $[n]^d$. They considered sets $A$ whose elements are chosen independently at random with probability $p$, and asked for which values of $p$ percolation is likely to occur. More precisely, let $P(G, r, p)$ denote the probability that $A$ percolates if $A$ is chosen with this distribution, and, noting that $P(G, r, p)$ is strictly increasing in $p$, define, for each $\alpha \in [0, 1]$,

$$p_\alpha = p_\alpha(G, r) := \inf\{p : P(G, r, p) \geq \alpha\}.$$

As is customary, we shall write $p_c$ for $p_{1/2}$, and call it the *critical probability*. Aizenman and Lebowitz [2] (see also Balogh and Pete [9]) showed that, for fixed $d$,

$$p_c([n]^d, 2) = \Theta\left(\frac{1}{\log n}\right)^{d-1},$$

and thus determined $p_c$ up to a constant factor when $r = 2$. Moreover, they showed that $p_\alpha$ satisfies the same relation for every fixed $\alpha \in (0, 1)$. When $r \geq 3$ the problem is somewhat harder, and it was not until 1999 that Cerf and Cirillo [12] determined the order of magnitude of $p_c([n]^3, 3)$. This result was later extended to all fixed $2 \leq r \leq d$ by Cerf and Manzo [13], who proved that

$$p_c([n]^d, r) = \Theta\left(\frac{1}{\log_{(r-1)} n}\right)^{d-r+1},$$

where $\log_{(r)}$ denotes an $r$-times iterated logarithm, $\log_{(r+1)}(n) = \log(\log_{(r)}(n))$. The bootstrap process has also been studied on the hypercube [3, 5, 6], on infinite trees [8] and on the random regular graph [10, 21], and has found applications in other areas: for example, techniques from [2], and more recently [6], were used to study the Ising model at zero temperature (see [15] and [22]).

Despite this extensive body of work, the threshold $p_c([n]^d, r)$ is known asymptotically only in the simplest case, $d = r = 2$. This important breakthrough was made by Holroyd [19], who proved that

$$p_c([n]^2, 2) = \frac{\pi^2}{18 \log n} + o\left(\frac{1}{\log n}\right).$$

(Here, and throughout, log refers to the natural logarithm, unless otherwise stated.) We shall discuss Holroyd's ideas in more detail later in the paper. Holroyd [20] also determined the critical constant in $d$ dimensions (for any fixed $d$) in the simpler "modified" bootstrap model, but was unable to do so for the standard case. Balogh and Bollobás [4] proved, using a general method of Friedgut and Kalai [17], that the threshold $p_c([n]^d, 2)$ undergoes a sharp transition in a weaker sense. More precisely, they showed that the critical window $p_{1-\varepsilon} - p_\varepsilon$ has width $o(p_c)$, but not that $(\log n)^{d-1} p_c$ converges. For $r \geq 3$, however, even this weaker result is unknown.

We shall determine the critical probability for 3-neighbor bootstrap percolation on $[n]^3$, up to a factor of $1 + o(1)$. In order to state our main result, we first need to define some functions. For each $k \in \mathbb{N}$, let

$$(1) \qquad \beta_k(u) := \frac{1}{2} - \frac{(1-u)^k}{2} + \frac{1}{2}\sqrt{1 + (4u - 2)(1-u)^k + (1-u)^{2k}},$$

so $\beta_k(u)^2 = (1 - (1-u)^k)\beta_k(u) + u(1-u)^k$, and let

$$(2) \qquad g_k(z) := -\log(\beta_k(1 - e^{-z})).$$

Now, for each $2 \leq r \leq d \in \mathbb{N}$, let

$$(3) \qquad \lambda(d, r) := \int_0^\infty g_{r-1}(z^{d-r+1}) \, dz.$$

The following theorem is the main result of this paper.

**Theorem 1.** *Let $\lambda(3, 3) \approx 0.4039$ be as defined above. Then*

$$p_c([n]^3, 3) = \frac{\lambda(3, 3) + o(1)}{\log \log n}$$

*as $n \to \infty$.*

The functions $\beta_1$, $g_1$ and $\lambda(2, 2)$ were introduced by Holroyd [19], who also showed that $\lambda(2, 2) = \pi^2/18$. We make the following conjecture, which is proved by Balogh et al. [7] in a forthcoming article.

**Conjecture 1.** *Let $d, r \in \mathbb{N}$, with $d \geq r \geq 2$. Then*

$$p_c([n]^d, r) = \left( \frac{\lambda(d, r) + o(1)}{\log_{(r-1)} n} \right)^{d-r+1}$$

*as $n \to \infty$.*

We remark (see Proposition 4 below) that $\lambda(d, r) < \infty$ for every $2 \leq r \leq d$, that $\lambda(3, 3) \approx 0.4039$ (by computer approximation), and that $d\lambda(d, d) \to$ $\pi^2/6$ as $d \to \infty$. Table 1 lists some approximate values of $\lambda(d, r)$ for $2 \le r \le 7$.

TABLE 1
*Values of $\lambda(d, r)$*

| | | | $d$ | | | |
|---|---|---|---|---|---|---|
| $r$ | **2** | **3** | **4** | **5** | **6** | **7** |
| 2 | 0.5483 | 0.9924 | 1.4797 | 1.9764 | 2.4760 | 2.9768 |
| 3 | – | 0.4039 | 0.8810 | 1.3864 | 1.8961 | 2.4078 |
| 4 | – | – | 0.3198 | 0.8024 | 1.3162 | 1.8338 |
| 5 | – | – | – | 0.2650 | 0.7431 | 1.2606 |
| 6 | – | – | – | – | 0.2265 | 0.6963 |
| 7 | – | – | – | – | – | 0.1979 |

The proof of Theorem 1 uses the techniques introduced by Cerf and Cirillo [12], Cerf and Manzo [13], and Holroyd [19], together with some new ideas. In particular, we shall need to introduce the following more general family of bootstrap processes.

Define a *bootstrap structure* $B(G, r(v))$ to be a graph $G$ together with a threshold function $r : V(G) \to \mathbb{N}$. Bootstrap percolation on such a structure is then defined in the obvious way, by setting $A_0 = A$ and

$$A_{t+1} := A_t \cup \{v \in V(G) : |\Gamma(v) \cap A_t| \ge r(v)\}$$

for each $t \ge 0$. This definition clearly includes all of the processes considered above; in particular, $B([n]^d, r)$ is the usual $r$-neighbor structure on the graph $[n]^d$.

We shall consider, in particular, the following two families of bootstrap structures. The first, which we shall call $C^*(n, 2)$, is an $n \times n \times 2$ cuboid, with threshold 2 in the "top" layer and threshold 3 in the "bottom" layer, that is, $r(v) = 2$ if $v \in [n]^2 \times \{1\}$ and $r(v) = 3$ if $v \in [n]^2 \times \{2\}$. The second, which we shall call $C(n, k)$, where $2 \le k \in \mathbb{N}$, is an $n \times n \times k$ cuboid, with threshold 2 in the top and bottom layers and threshold 3 in each of the $k - 2$ middle layers. Both $C^*(n, 2)$ and $C(n, k)$ have the edges induced by the lattice $\mathbb{Z}^3$.

We need two more definitions. Say that $A \subset C^*(n, 2)$ *semi-percolates* if $[A]$ contains all vertices with threshold 2, and, for each $\alpha \in (0, 1)$, write $p_\alpha^{(s)}(C^*(n, 2)) := \inf\{p : \mathbb{P}(A \text{ semi-percolates}) \ge \alpha\}$. Recall also from above that $p_\alpha(C(n, k))$ is defined similarly for percolation (i.e., full occupation).

We shall prove Theorem 1 using the following result.

THEOREM 2. *For every $\varepsilon > 0$, there exists $K = K(\varepsilon) \in \mathbb{N}$ such that, if $k \ge K$, $\delta = \delta(k) > 0$ is sufficiently small, and $n = n(k, \delta, \varepsilon) \in \mathbb{N}$ is sufficiently*

*large, then*

$$\frac{\lambda(3,3) - \varepsilon}{\log n} \le p_\delta(C(n,k)) \le p_{1-\delta}^{(s)}(C^*(n,2)) \le \frac{\lambda(3,3) + \varepsilon}{\log n}.$$

In fact, we shall need a somewhat more technical statement in order to deduce Theorem 1 (see Corollary 14 below), but Theorem 2 is morally what is required.

Being a little imprecise, the basic idea of the proof of Theorem 1 is as follows. With high probability, complete occupation occurs in bootstrap percolation on $B([n]^3, 3)$ (roughly) if and only if there exists, somewhere in $[n]^3$, a cuboid $R$, with side-lengths about $\log n$, which is internally spanned, that is, $R$ satisfies $[A \cap R] = R$. If this exists, we (rather vaguely) refer to such a cuboid $R$ as a "critical droplet." We couple the process occurring on the sides (i.e., faces) of the droplet $R$ in two different ways, in order to prove upper and lower bounds on the probability that $R$ "grows sideways." To prove the upper bound in Theorem 1, we couple with $C^*(\log n, 2)$, and use the upper bound in Theorem 2; for the lower bound we couple with $C(\log n, k)$ and use a counting argument as in [12, 13]. This allows us to show that no large connected component of infected sites forms anywhere in $[n]^3$.

Several of our lemmas generalize easily to a multidimensional setting, and will be used (in this more general form) in [7]. We shall therefore often work in $[n]^d \times [k]^\ell$ for general $d \ge 2$ and $\ell \ge 0$. However, the reader should always be thinking of the case $d = 2$ and $\ell = 1$, and our terminology will reflect this.

The paper is organized as follows: In Section 2 we collect some of the definitions and basic tools which we shall use throughout the paper, and prove some simple bounds on $\lambda(d, r)$. In Section 3 we prove the upper bound in Conjecture 1 (and hence in Theorem 1 also), and in Section 4 we prove Theorem 2 (and Corollary 14). Finally, in Section 5, we deduce the lower bound in Theorem 1.

## 2. Tools and notation.

In this section we shall make various definitions and introduce some notation which we shall use throughout the paper. We have labeled some of these in order to highlight those that are most crucial.

We begin by defining some slightly more general versions of the bootstrap structures described above, which we shall call $C([n]^d \times [k]^\ell, r)$ and $C^*([n]^d \times [2]^\ell, r)$. We think of $[n]^d \times [k]^\ell$ as a box $[n]^d$ of "thickness" $[k]^\ell$.

DEFINITION. Let $n, d, \ell, r \in \mathbb{N}_0$, with $2 \le r \le d$. Then $C^*([n]^d \times [2]^\ell, r)$ is the bootstrap structure such that:

(a) the vertex set is $[n]^d \times [2]^\ell$,
(b) the edge set is induced by $\mathbb{Z}^{d+\ell}$,
(c) $v = (a_1, \ldots, a_d, b_1, \ldots, b_\ell)$ has threshold $r$ if $b_j = 1$ for each $j \in [\ell]$,

(d)  $v = (a_1, \ldots, a_d, b_1, \ldots, b_\ell)$ has threshold $r + \ell$ otherwise.

We say that a set $A \subset [n]^d \times [2]^\ell$ *semi-percolates* in $C^*([n]^d \times [2]^\ell, r)$ if $[A]$ contains all vertices with threshold $r$. Note that $C^*(n, 2) = C^*([n]^2 \times [2], 2)$.

DEFINITION.    Let $n, d, k, \ell, r \in \mathbb{N}_0$, with $2 \le r \le d$ and $k \ge 2$. Then $C([n]^d \times [k]^\ell, r)$ is the bootstrap structure such that:

(a)  the vertex set is $[n]^d \times [k]^\ell$,
(b)  the edge set is induced by $\mathbb{Z}^{d+\ell}$,
(c)  $v = (a_1, \ldots, a_d, b_1, \ldots, b_\ell)$ has threshold $r + |\{j \in [\ell] : b_j \notin \{1, k\}\}|$.

Note that $C(n, k) = C([n]^2 \times [k], 2)$ for any $k \ge 2$, and $B([n]^d, r) = C([n]^d \times [k]^0, r)$.

We next, with the case $d = 2$ in mind, define a *rectangle* $R$ in $[n]^d \times [k]^\ell$ to be a set

$$[(a_1, \ldots, a_d), (b_1, \ldots, b_d)] := \{(x_1, \ldots, x_d, y_1, \ldots, y_\ell) : x_i \in [a_i, b_i], y_i \in [k]\}.$$

We also identify these with rectangles in $[n]^d = [n]^d \times [k]^0$ in the obvious way. The *dimensions* of $R$ is the vector

$$\dim(R) := (b_1 - a_1 + 1, \ldots, b_d - a_d + 1) \in \mathbb{N}^d$$

and the *semi-perimeter* of $R$ is

$$\phi(R) := \sum_i (b_i - a_i + 1).$$

The longest side-length of $R$ is $\mathrm{long}(R) := \max\{b_i - a_i + 1\}$, and the shortest side-length of $R$ is $\mathrm{short}(R) := \min\{b_i - a_i + 1\}$.

A *component* of a set $S \subset \mathbb{Z}^d$ is a maximal connected set in the graph $\mathbb{Z}^d[S]$ (the subgraph of $\mathbb{Z}^d$ induced by $S$), and the *diameter* of $S$ is

$$\mathrm{diam}(S) := \sup_{x,y} \{\|x - y\|_\infty + 1 : x \text{ and } y \text{ are in the same component of } S\}.$$

Note that if $S$ is a rectangle in $[n]^d \times [k]^\ell$, then $\mathrm{diam}(S) = \max\{\mathrm{long}(R), k\}$.

Let $S \subset [n]^d \times [k]^\ell$. The *projection* $\Pi(S) \subset [n]^d$ of $S$ is the set

$$\Pi(S) := \{\mathbf{x} \in [n]^d : (\mathbf{x}, \mathbf{y}) \in S \text{ for some } \mathbf{y} \in [k]^\ell\}.$$

We have defined the completion $[A]$ of $A$; now we shall define the *span*, $\langle A \rangle$. We emphasize that this is *not* the usual definition. First, note that for each subset $S \subset [n]^d$ there is a smallest rectangle, $R(S)$, such that $S \subset R(S)$.

DEFINITION. Let $n, k \in \mathbb{N}$ and $A \subset C([n]^d \times [k]^\ell, r)$. Let $C_1, \ldots, C_m$ denote the collection of connected components in $\Pi([A])$. The *span of $A$* is defined to be the following collection of rectangles:

$$\langle A \rangle := \{R(C_1), \ldots, R(C_m)\}.$$

If $\Pi([A])$ is connected (i.e., $m = 1$), then *$A$ spans the rectangle $R(C_1)$*. Also, if $\langle A' \rangle = \{R\}$ for some $A' \subset A$, then *$A$ internally spans $R$*.

If $\langle A \rangle = \{R\}$, that is, $A$ spans $R$, then we shall usually write simply $\langle A \rangle = R$. Note that $A \subset R$ internally spans $R$ if and only if $R \in \langle A \rangle$.

We now prove some simple properties of the functions $\beta_k(u)$, $g_k(u)$ and $\lambda(d, r)$ defined in Section 1.

PROPOSITION 3. *Let $k \in \mathbb{N}_0$, and $\beta_k$ and $g_k$ be the functions defined in (1) and (2):*

(a) *$\beta_k(u)$ is increasing in $u$ on $[0, 1]$ and $g_k(u)$ is decreasing in $u$ on $(0, \infty)$.*

(b) *$\beta_{k+1}(u) \geq \beta_k(u) \in [0, 1]$ for $u \in [0, 1]$, and $g_{k+1}(u) \leq g_k(u)$ for $u \in (0, \infty)$.*

(c) *$g_k(z) \leq 2e^{-zk}$ if $z$ is sufficiently large.*

PROOF. We use Lemma 6 below, which says that the probability $L_k(m, u)$ that there is no "$L$-gap" (defined below) in a sequence of events of length $m$, with each event having probability $u \in [0, 1]$, satisfies

$$\beta_{k+1}(u)^{m+1} \leq L_k(m, u) \leq \beta_{k+1}(u)^m.$$

(Note that the proof of Lemma 6 is straightforward and self-contained.) It is clear from the definition of $L$-gaps that $L_k(m, u)$ is strictly increasing in both $k$ and $u$. Thus, applying the displayed equation for sufficiently large $m$, it follows that $\beta_k$ is increasing in both $u$ and $k$, and that $\beta_k(u) \in [0, 1]$. The facts about $g_k$ in parts (a) and (b) now follow from those about $\beta_k$ by (2), the definition of $g_k$.

For part (c) recall that $\beta_k(u)^2 = (1 - (1-u)^k)\beta_k(u) + u(1-u)^k$, so $\beta_k(u) \geq 1 - (1-u)^k$ for $u \in [0, 1]$. Recall also that $-\log(1-x) \leq 2x$ if $x > 0$ is sufficiently small. Therefore,

$$-\log(\beta_k(1 - e^{-z})) \leq -\log(1 - e^{-zk}) \leq 2e^{-zk},$$

if $z$ is sufficiently large, as required. $\square$

Although we cannot solve the integral (3) exactly, the following proposition gives some bounds on $\lambda(d, r)$. The proofs are all straightforward, so we give only a sketch.

PROPOSITION 4.    *Let $2 \le r \le d$ and let $\lambda(d, r)$ be the function defined in* [(3)](3):

(a) $\lambda(d, r) < \infty$.
(b) $\lambda(d, 2) = \frac{d-1}{2} + o(1)$.
(c) $d\lambda(d, d) \to \frac{\pi^2}{6}$ *as* $d \to \infty$.

PROOF.    First we show that, for every $d, k \in \mathbb{N}$,

$$\int_0^\infty g_k(z^d)\, dz \le \int_0^\infty g_1(z^d)\, dz = \frac{d}{2} + o(1).$$

The first inequality follows by Proposition [3](3)(b). For the second we use the fact that $\beta_1$ is increasing and continuous, and the following simple facts:

(i) $(1 - e^{-x})/x \to 1$ as $x \to 0$.
(ii) $\beta_1(x)/\sqrt{x} \to 1$ as $x \to 0$.
(iii) $\int_1^\infty g_1(z^d)\, dz \to 0$ as $d \to \infty$.

It follows that, letting $d \to \infty$,

$$\beta_1(1 - e^{-z^d}) = (1 + o(1))\beta_1(z^d) = (1 + o(1))\sqrt{z^d}$$

for every $z \in (0, 1)$, and hence,

$$\int_0^\infty g_1(z^d)\, dz = \int_0^1 \left( -\frac{d}{2} \log z + o(1) \right) dz + \int_1^\infty g_1(z^d)\, dz = \frac{d}{2} + o(1),$$

as required. Parts (a) and (b) now follow by the definition of $\lambda(d, r)$.

For part (c), note that

$$\beta_k(1 - e^{-z}) = \frac{1}{2} - \frac{e^{-zk}}{2} + \frac{1}{2}\sqrt{(1 - e^{-zk})^2 + 4e^{-zk}(1 - e^{-z})},$$

and that $x \le \sqrt{x^2 + y} \le x + \frac{y}{2x}$ if $y \ge 0$. Thus, since $1 - e^{-z} \le z$,

$$1 - e^{-zk} \le \beta_k(1 - e^{-z}) \le 1 - e^{-zk} + \frac{ze^{-zk}}{1 - e^{-zk}}.$$

Hence, making the substitution $x = e^{-zk}$,

$$k\int_0^\infty g_k(z)\, dz \ge \int_0^\infty -k\log(1 - e^{-zk})\, dz = \int_0^1 -\log(1 - x)\frac{dx}{x} = \frac{\pi^2}{6}$$

by [18], number 4.291.2. Moreover, $\frac{z}{1 - e^{-zk}}$ is increasing on $z > 0$ (by simple calculus). Thus, letting $k$ be large, $\delta := \frac{1}{\sqrt{k}}$, and using the substitution $x = (1 - 2\delta)e^{-zk}$,

$$k\int_0^\infty g_k(z)\, dz \ge \int_{\delta/k}^\delta -k\log\left( 1 - e^{-zk} + \frac{\delta e^{-zk}}{1 - e^{-\delta k}} \right) dz$$

$$\geq \int_{\delta/k}^{\delta} -k \log(1 - e^{-zk} + 2\delta e^{-zk})\, dz$$

$$\geq \int_{e^{-\delta k}}^{(1-2\delta)e^{-\delta}} -\log(1-x)\frac{dx}{x} \to \frac{\pi^2}{6}$$

as $k \to \infty$, as required.  □

REMARK 1. The constant determined by Holroyd [20] for the modified bootstrap model is $\pi^2/6$. Thus, when $d$ is large, $\lambda(d,d)$ differs from the critical constant in that model by a factor of $d$. It is tempting to suggest a simple explanation for this: a blocking set (i.e., an $L$-gap, see Section 3) is $d$ times larger in bootstrap percolation than in modified bootstrap (see Lemmas 6 and 7 below). Caution is required, however, since this heuristic fails when the number of dimensions is small; in particular, when $d = 2$ the critical constants differ instead by a factor of three.

We next state the FKG (Fortuin–Kasteleyn–Ginibre) inequality [16], the van den Berg–Kesten Lemma [11] and Reimer's Theorem [23]. We shall use the former two results several times; the latter will only be used once, but will play a key role in the proof. We remark that although Reimer's Theorem appears rather naturally in our proof, it can be avoided; indeed, the proof in [7] of the (more general) Conjecture 1 does not use it.

We begin with the simplest of the three results, the FKG inequality. Let $E \colon \mathcal{P}(n) \to \{T, F\}$ be an *event* defined on the cube $\mathcal{P}(n) = \{0,1\}^n$, that is, $E$ is a subset of $\{0,1\}^n$. $E$ is said to be *increasing* if, for any two sets $X, Y \subset [n]$,

$$E(X) \wedge (X \subset Y) \Rightarrow E(Y).$$

We write $\mathbb{P}_p$ for the product measure on $\{0,1\}^n$ with $\mathbb{P}_p(j \in A) = p$ for each $j \in [n]$.

THE FKG INEQUALITY. *Let $n \in \mathbb{N}$ and $p \in (0,1)$, and let $E$ and $F$ be increasing events on the cube $\{0,1\}^n$. Then*

$$\mathbb{P}_p(E \cap F) \geq \mathbb{P}_p(E)\mathbb{P}_p(F).$$

Now let $E$ and $F$ be two events defined on $\mathcal{P}(n)$, and let $S \subset [n]$. A *witness set* for the event "$E(S)$ holds" is a disjoint pair of sets $(U, V)$ such that $U \subset S$, $S \cap V = \varnothing$, and

$$(U \subset X) \wedge (X \cap V = \varnothing) \Rightarrow E(X)$$

for any set $X \in \mathcal{P}(n)$. The events $E$ and $F$ are said to *occur disjointly* at a point $S \in \mathcal{P}(n)$ if there exist witness sets $(U, V)$ and $(U', V')$ for the events

"$E(S)$ holds" and "$F(S)$ holds" respectively, such that the sets $U \cup V$ and $U' \cup V'$ are disjoint.

We write $E \circ F$ for the event that $E$ and $F$ occur disjointly. The following lemma is an important and much-used tool in percolation theory, and was proved by van den Berg and Kesten [11].

VAN DEN BERG–KESTEN LEMMA. *Let $n \in \mathbb{N}$ and $p \in (0, 1)$, and let $E$ and $F$ be increasing events defined on the cube $\{0, 1\}^n$. Then*

$$\mathbb{P}_p(E \circ F) \le \mathbb{P}_p(E)\mathbb{P}_p(F).$$

The following substantial generalization of the van den Berg–Kesten Lemma was conjectured by van den Berg and Kesten [11] and proved by Reimer [23].

REIMER'S THEOREM. *Let $n \in \mathbb{N}$ and $p \in (0, 1)$, and let $E$ and $F$ be arbitrary events defined on the cube $\{0, 1\}^n$. Then*

$$\mathbb{P}_p(E \circ F) \le \mathbb{P}_p(E)\mathbb{P}_p(F).$$

We conclude the section with a little more notation. Given a set $S$, and $p \in [0, 1]$, say that $A \in \text{Bin}(S, p)$ if the elements of $A \subset S$ are chosen independently at random with probability $p$. If $R$ is a rectangle in $C([n]^d \times [k]^\ell, r)$, then let

$$P_p(R) := \mathbb{P}(R \in \langle A \rangle \mid A \in \text{Bin}(R, p)),$$

that is, the probability that $A \in \text{Bin}(R, p)$ spans $R$.

A set is said to be *occupied* if it is nonempty (i.e., contains some element of $A$), and it is said to be *full* if every site is in $A$. We shall use throughout the paper the notation

$$q := -\log(1 - p)$$

as in [19]. Note that $p \sim q$ for small $p$. The advantage of this notation is the fact that

$$(4) \qquad\qquad \beta_k(1 - (1 - p)^n) = e^{-g_k(nq)}.$$

For any $a, b \in \mathbb{Z}$, we write $[a, b]$ for the set $\{n \in \mathbb{Z} : a \le n \le b\}$, and $[a] = [1, a]$. Given two functions $f, g \colon \mathbb{N} \to \mathbb{R}$, we say that $f \gg g$ if $f(n)/g(n) \to \infty$ as $n \to \infty$. Whenever $v$ is a vector, $v_j$ will be its $j$th coordinate. Finally, if $G$ is an oriented tree, then $\vec{\Gamma}(u) = \vec{\Gamma}_G(u) := \{v \in V(G) \colon u \to v\}$.

**3. A general upper bound.** We begin by proving the upper bound in Theorem 1. The proof is straightforward (though slightly technical); since essentially the same method gives the upper bound in Conjecture 1 for all $d$ and $r$, we shall give the general argument. We refer the reader to [19] and [20] (see also [2, 27] and [25]), where many of the ideas we shall use originated.

THEOREM 5. *Let $d, \ell, r \in \mathbb{Z}$, with $d \geq r \geq 2$ and $\ell \geq 0$, and let $\varepsilon > 0$. Suppose $n \in \mathbb{N}$,*

$$p \geq \left( \frac{\lambda(d + \ell, \ell + r) + \varepsilon}{\log_{(r-1)} n} \right)^{d-r+1},$$

*and the elements of $A \subset C^*([n]^d \times [2]^\ell, r)$ are chosen independently at random with probability $p$. Then*

$$\mathbb{P}(A \text{ semi-percolates in } C^*([n]^d \times [2]^\ell, r)) \to 1$$

*as $n \to \infty$. In particular,*

$$p_c([n]^d, r) \leq \left( \frac{\lambda(d, r) + o(1)}{\log_{(r-1)} n} \right)^{d-r+1}.$$

The idea of the proof is quite simple, but the details are a little technical, so we begin with some motivation. It follows from [13] that a "critical droplet" in $B([n]^d, r)$ has size roughly $[\log n]^d$. Suppose a cube $R$ of this size is completely infected, and consider two "layers" next to it in some direction, which form a copy $S$ of $[\log n]^{d-1} \times [2]$. Since all vertices of $R$ are infected, each vertex in the layer of $S$ adjacent to $R$ has one already-infected neighbor, and so requires only $r - 1$ more infected neighbors from inside $S$. Thus, if the set $A \cap S$ semi-percolates (i.e., completely occupies the layer with threshold $r - 1$) in the bootstrap structure $C^*([\log n]^{d-1} \times [2], r - 1)$, then the cube $R$ will grow sideways by one step. Hence, the critical droplet is likely to grow if percolation is likely to occur in $C^*([\log n]^{d-1} \times [2], r - 1)$.

Applying the same logic to the structure $C^*([\log n]^{d-1} \times [2], r - 1)$, we see that a critical droplet has size about $[\log \log n]^{d-1} \times [2]$. (It is important to note that we only require the "top layer," i.e., the vertices with threshold $r - 1$, to be infected.) This droplet is likely to grow if semi-percolation is likely to occur in the copies of $C^*([\log \log n]^{d-2} \times [2]^2, r - 2)$ on its sides, and so on. Iterating $r - 2$ times, we see that percolation is likely to occur in $B([n]^d, r)$ if semi-percolation is likely to occur in $C^*([\log_{(r-1)} n]^{d-r+2} \times [2]^{r-2}, 2)$.

We begin with the base case, $C^*([n]^d \times [2]^\ell, 2)$. Our first aim is to give a lower bound on the probability that the set $A$ semi-percolates in $C^*([n]^d \times [2]^\ell, 2)$ by considering one particular way in which the percolation may occur.

Slightly more precisely, we shall consider the growth of an infected "droplet" which begins life in the bottom left-hand corner, and grows upward and rightward by "crossing rectangles" as follows.

Let $e_t$ be the site $(0, \ldots, 0, 1, 0, \ldots, 0) \in \mathbb{Z}^{d+\ell}$ with a single 1 in position $t$. We shall sometimes write $\mathbf{1}^\ell$ to denote the vector $(1, \ldots, 1) \in [2]^\ell$.

DEFINITION. Let $2 \leq n, d \in \mathbb{N}$, $\ell \in \mathbb{N}_0$, $R \subset C^*([n]^d \times [2]^\ell, 2)$ be a rectangle, and $A \subset [n]^d \times [2]^\ell$. For each $t \in [d]$, let $R_t^+ := \{v \notin R : v - e_t \in R, r(v) = 2\}$ be the set of vertices with threshold 2 immediately to the right of $R$, and $R_t^- := \{v \notin R : v + e_t \in R, r(v) = 2\}$ be those immediately to the left of $R$, both right and left being in direction $t$.

Now, let

$$A_t^R := (A \cap (R \cup R_t^+)) \cup R_t^- .$$

We say $R$ is *semi-crossed from left to right in direction $t$* by $A$ if the set $[A_t^R]$ contains all vertices in $R$ with threshold 2.

In other words, if $R$ is semi-crossed by $A$, and the sites with threshold 2 to the left of $R$ have already been infected, the sites of $R$ with threshold 2 will then also be infected.

In order to bound the probability that a rectangle in $C^*([n]^d \times [2]^\ell, 2)$ is semi-crossed, we need to introduce the concept of an $L$-gap in a sequence of events.[4] Let $\ell, m \in \mathbb{N}$, and consider some sequence of events

$$\mathcal{E} = \{U_i : i \in [m+1]\} \cup \{V_j^{(i)} : i \in [\ell], j \in [m]\}.$$

An $L$-gap in $\mathcal{E}$ is an event $\neg(U_i \vee U_{i+1} \vee V_i^{(1)} \vee \cdots \vee V_i^{(\ell)})$ for some $i \in [m]$.

LEMMA 6. *Let $\ell, m \in \mathbb{N}$, let $u \in (0, 1)$, and suppose that each event in the set*

$$\mathcal{E} = \{U_i : i \in [m+1]\} \cup \{V_j^{(i)} : i \in [\ell], j \in [m]\}$$

*occurs independently with probability $u$.*

*Let $L(m, u)$ denote the probability that there is no $L$-gap in $\mathcal{E}$. Then*

$$\beta_{\ell+1}(u)^{m+1} \leq L(m, u) \leq \beta_{\ell+1}(u)^m,$$

*where $\beta_{\ell+1}(u)$ is the function defined in the Introduction.*

PROOF. We partition the event that there is no $L$-gap in $\mathcal{E}$ into three cases, and use induction on $m$. Let $L(-1, u) = L(0, u) = 1$, and note that $\beta_{\ell+1}(u) \in (0, 1)$ for $u \in (0, 1)$, so the induction hypothesis holds for $m \in$ $\{-1, 0\}$. So let $m \in \mathbb{N}$, and observe that either at least one of the events $U_1, V_1^{(1)}, \ldots, V_1^{(\ell)}$ occurs, or none of these occurs but $U_2$ does, or none does and $U_2$ also does not. Conditional on these events, the probabilities that there is no $L$-gap are $L(m-1, u)$, $L(m-2, u)$ and $0$, respectively. Thus,

$$L(m, u) = (1 - (1 - u)^{\ell+1})L(m-1, u) + u(1 - u)^{\ell+1}L(m-2, u)$$

for every $m \geq 1$. Furthermore,

$$\beta_{\ell+1}(u)^2 = (1 - (1 - u)^{\ell+1})\beta_{\ell+1}(u) + u(1 - u)^{\ell+1},$$

and so the result follows by induction, as claimed.  □

---

[4]Note that $L$-gaps are so named because of their shape; this $L$ is not a variable.

We now deduce the following bound on the probability that a rectangle is semi-crossed.

LEMMA 7.  *Let* $2 \leq n, d \in \mathbb{N}$, $\ell \in \mathbb{N}_0$, $R \subset C^*([n]^d \times [2]^\ell, 2)$ *be a rectangle,* $p \in (0, 1)$ *and* $A \in \mathrm{Bin}(R, p)$. *Let* $t \in [d]$ *and write* $v(t) = \prod_{i \neq t} a_i$ *and* $u(t) = 1 - (1 - p)^{v(t)}$, *where* $\dim(R) = (a_1, \ldots, a_d)$. *Then*

$$\mathbb{P}(R \text{ is semi-crossed in direction } t \text{ by } A) \geq \beta_{\ell+1}(u(t))^{a_t+1}.$$

PROOF.  Consider the sequence of events

$$\mathcal{E} = \{U_i : i \in [a_t + 1]\} \cup \{V_j^{(i)} : i \in [\ell], j \in [a_t]\},$$

where

$$U_i = \{[a_1] \times \cdots \times [a_{t-1}] \times \{i\} \times [a_{t+1}] \times \cdots \times [a_d] \times \mathbf{1}^\ell \text{ is occupied}\}$$

and

$$V_j^{(i)} = \{[a_1] \times \cdots \times [a_{t-1}] \times \{j\} \times [a_{t+1}] \times \cdots \times [a_d] \times (\mathbf{1}^\ell + e_i) \text{ is occupied}\}.$$

As before, a set is said to be occupied if it contains at least one element of the set $A \subset R$. Note that each of the events $U_i$ and $V_j^{(i)}$ occurs independently with probability $u(t)$.

We claim that if $\mathcal{E}$ has no $L$-gap, then $R$ is crossed from left to right in direction $t$ by $A$. Indeed, let $m$ be the minimal index such that some element of $\{v \in R : v_t = m, r(v) = 2\}$ is not in $[A_t^R]$. Then the event

$$\neg(U_m \vee U_{m+1} \vee V_m^{(1)} \vee \cdots \vee V_m^{(\ell)})$$

holds, and is an $L$-gap. The result now follows by Lemma 6.  □

For each $p \in (0, 1)$, and each $n, d, \ell, r \in \mathbb{N}_0$ with $2 \leq r \leq d$, let $P(n, d, \ell, r, p)$ denote the probability that a set $A \in \mathrm{Bin}(C^*([n+1]^d \times [2]^\ell, r), p)$ semi-percolates in $[n]^d \times [2]^\ell$, that is, $[n]^d \times \mathbf{1}^\ell \subset [A]$. Note that we are allowed to use active sites in the layer outside $[n]^d \times [2]^\ell$; this technicality will help to simplify the proof below.

The next lemma gives a lower bound on $P(n, d, \ell, 2, p)$ by considering one way in which the spanning could occur. Let

$$(5) \qquad G(d, \ell, 2, p) := \frac{d \int_0^\infty g_{\ell+1}(z^{d-1}) \, dz}{p^{1/(d-1)}} = \frac{d\lambda(d+\ell, \ell+2)}{p^{1/(d-1)}}.$$

LEMMA 8.  *Let* $2 \leq n, d \in \mathbb{N}$, $\ell \in \mathbb{N}_0$, $\varepsilon > 0$, *and* $p > 0$ *be sufficiently small.* *Then*

$$P(n, d, \ell, 2, p) \geq \exp(-(1+\varepsilon)G(d, \ell, 2, p)).$$

REMARK 2.  Note that $n$ does not appear in the expression on the right-hand side.

PROOF OF LEMMA 8.  Let $n$, $d$, $\ell$, $\varepsilon$ and $p$ be as described, and let $A \in \text{Bin}(C^*([n+1]^d \times [2]^\ell, 2), p)$. We shall describe sufficient conditions for semi-percolation to occur. Let $r := \lfloor p^{1/(2d-2)} \rfloor$, and let $E_0$ denote the event that the set

$$M_0 := \bigcup_{j=1}^d \{\mathbf{x} \in [n]^d \times \mathbf{1}^\ell : x_j \in [r] \text{ and } x_i = 1 \text{ if } i \neq j\}$$

is full, that is, $M_0 \subset A$. Now, for each $j \in [d]$ and $t \in \mathbb{N}$, let

$$R_j(t) := \{\mathbf{x} \in [n]^d : x_j \in [tr+1, (t+1)r] \text{ and } x_i \in [tr] \text{ if } i \neq j\} \times [2]^\ell,$$

and note that $R_j(t)$ is a $[tr]^{d-1} \times [r] \times [2]^\ell$-cuboid. Let $E_j(t)$ denote the event that $R_j(t)$ is left-to-right semi-crossed in direction $j$ by $A$.

We first claim that if $E_0$ and $E_j(t)$ hold for every $j \in [d]$ and $1 \leq t \leq n/r$, then $A$ semi-percolates in $C^*([n]^d \times [2]^\ell, 2)$. Indeed, $[r]^d \times \mathbf{1}^\ell \subset [A]$ since $E_0$ holds, so $R_j(1) \subset [A]$ for each $j \in [d]$ since $E_j(1)$ holds. But if $R_j(1) \subset [A]$ for each $j \in [d]$, then $[2r]^d \times \mathbf{1}^\ell \subset [A]$. Repeating this argument shows that $[tr]^d \times \mathbf{1}^\ell \subset [A]$ for each $t \leq n/r$.

It remains to bound the probability of the events $E_0$ and $E_j(t)$; since these events are all increasing, by the FKG inequality we may bound the probability of their intersection from below by the product of their probabilities. [Note that they are not independent, since the event $E_j(t)$ depends on the set $A \cap (R_j(t) \cup (R_j(t))_j^+)$.] It is easy to see that $\mathbb{P}(E_0) \geq p^{dr}$, and by Lemma 7 and (4) we have

$$\mathbb{P}(E_j(t)) \geq \beta_{\ell+1}(1 - (1-p)^{(tr)^{d-1}})^{r+1} = \exp(-(r+1)g_{\ell+1}(q(tr)^{d-1})).$$

Note also that

$$\sum_{t=1}^\infty g_{\ell+1}(q(tr)^{d-1}) \leq \frac{1}{q^{1/(d-1)}r} \int_0^\infty g_{\ell+1}(z^{d-1}) \, dz$$

since $g_{\ell+1}(z)$ is decreasing on $(0, \infty)$. Thus,

$$\mathbb{P}(A \text{ semi-percolates})$$

$$\geq \mathbb{P}(E_0) \prod_{j,t} \mathbb{P}(E_j(t))$$

$$\geq p^{dr} \exp\left(-d(r+1) \sum_{t=1}^{\infty} g_{\ell+1}(q(tr)^{d-1})\right)$$

$$\geq \exp\left(-\frac{d}{p^{1/(2d-2)}} \log\left(\frac{1}{p}\right) - \frac{d(r+1)}{q^{1/(d-1)}r} \int_0^{\infty} g_{\ell+1}(z^{d-1}) \, dz\right)$$

$$\geq \exp\left(-\frac{(1+\varepsilon)d\lambda(d+\ell,\ell+2)}{p^{1/(d-1)}}\right)$$

if $p$ is sufficiently small (as a function of $d$ and $\varepsilon$), as required. The last inequality holds by (3), and because $p \sim q$, $r+1 \sim r$ and $p^{-1/(2d-2)} \log(1/p) \ll p^{-1/(d-1)}$ as $p \to 0$. $\square$

Now we use Lemma 8 to prove Theorem 5 in the case $r = 2$.

LEMMA 9. *Let* $2 \leq d \in \mathbb{N}$, $\ell \in \mathbb{N}_0$, $\varepsilon > 0$, *and*

$$p = \left(\frac{\lambda(d+\ell,\ell+2) + \varepsilon}{\log n}\right)^{d-1}.$$

*Then*

$$P(n, d, \ell, 2, p) \to 1$$

*as* $n \to \infty$.

PROOF. Let $n$, $d$ and $\ell$ be as given, and assume $\varepsilon$ is sufficiently small. Recall that the function $G$ was defined in (5), and note that we have chosen $p$ so that $G(d, \ell, 2, p) \leq (1 - \varepsilon^2)d \log n$. Let $m = \lceil p^{-3/(d-1)} \rceil$, partition $[n]^d \times [2]^\ell$ into blocks of size $[m]^d \times [2]^\ell$, and run $\lfloor n/m \rfloor^d$ independent bootstrap processes on the intersection of $A$ with each block. By Lemma 8, the probability that $A$ semi-percolates in at least one of them is at least

$$1 - (1 - \exp(-(1+\varepsilon^3)G(d,\ell,2,p)))^{\lfloor n/m \rfloor^d} \geq 1 - (1 - n^{-(1-\varepsilon^3)d})^{\lfloor n/m \rfloor^d} \to 1$$

as $n \to \infty$, since $n^{\varepsilon^3} \gg m$.

Next, consider all $[m]^{d-1} \times [1]^{\ell+1}$ cuboids in $[n]^d \times [2]^\ell$. If $A$ semi-percolates in some $[m]^d \times [2]^\ell$ block, but $A$ does not semi-percolate in $C^*([n]^d \times [2]^\ell, 2)$, then one of these cuboids must be empty. But $pm^{d-1} \geq p^{-2}$, so the probability that at least one is empty is at most

$$dn^d(1-p)^{m^{d-1}} \leq dn^d \exp(-pm^{d-1}) \leq dn^d \exp\left(-\left(\frac{\log n}{\lambda + \varepsilon}\right)^2\right) \to 0$$

as $n \to \infty$, where $\lambda = \lambda(d + \ell, \ell + 2)$. Since the events "there exists an $[m]^d \times [2]^\ell$ block in which $A$ semi-percolates" and "all the $[m]^{d-1} \times [1]^{\ell+1}$ blocks in $[n]^d \times [2]^\ell$ are occupied" are both increasing events, the result follows by the FKG inequality. $\square$

Having proved the base case, the general result follows by a well-known and standard method (see [20], for example). We use the following two straightforward lemmas, which somewhat simplify the proof.

LEMMA 10 (Lemma 2 of [20]).    *For any $d \geq 3$, $\ell \geq 0$ and $\varepsilon > 0$, if $n$ is sufficiently large and $p^{-2d} \leq n^\varepsilon$, then*

$$P(n, d, \ell, 3, p) \geq \exp(-n^{1+\varepsilon}).$$

SKETCH OF PROOF.    Let $m = \frac{2d \log n}{p}$, and consider the set

$$M = \bigcup_{j=1}^{d} \{\mathbf{x} \in [n]^d \times \mathbf{1}^\ell : x_j \in [n] \text{ and } x_i \in [m] \text{ if } i \neq j\}.$$

The probability that $M \subset A$ is at least $p^{dm^{d-1}n} \geq \exp(-n^{1+\varepsilon})$, and the probability that $A$ semi-percolates given $M \subset A$ is at least the probability that every $[m] \times [1]^{d+\ell-1}$ cuboid is occupied, which is at least $1 - dn^d e^{-pm} = 1 - o(1)$. $\square$

The next lemma was proved in [2] for $G = [n]^d$ and $r = 2$, but the proof generalizes easily to our case.

LEMMA 11 (Aizenman and Lebowitz [2]).    *For each $2 \leq r \leq d \in \mathbb{N}$ and $\ell \in \mathbb{N}_0$, there exists $\delta = \delta(d, \ell, r) > 0$ and $C = C(d, \ell, r) < \infty$ such that, if $P(m, d, \ell, r, p) \geq 1 - \delta$, then*

$$P(n, d, \ell, r, p) \geq 1 - Ce^{-n/m}$$

*for every $n \geq m$.*

SKETCH OF PROOF.    First note that, by taking $C$ large, we may assume that $n/m$ is sufficiently large, since otherwise the result is trivial. The idea is to partition $[n]^d \times [2]^\ell$ into blocks of size (roughly) $[m]^d \times [2]^\ell$, and run the bootstrap process independently in each block. Call a block $B$ "active" if $A \cap B$ semi-percolates in $B$, and "inactive" otherwise.

Suppose that $A$ does not semi-percolate in $[n]^d$. We claim that every connected component of inactive $[m]^d$-blocks must span two opposite sides of $[n]^d$, that is, must touch both faces of $[n]^d$ in some direction. Indeed, consider a component $X$ of inactive $[m]^d$-blocks which does not span two opposite sides of $[n]^d$. Let $Y$ denote the collection of sites with threshold $r$ which are in blocks of $X$, but which are not in $[A]$. We claim that $Y$ is empty, and hence that $X$ is empty.

Note that if a block $B$ is on the boundary of $X$ (but not in $X$), then it is active, since $X$ is a component. Thus, if $x \in B$, and $r(x) = r$, then $x \in [A]$. Hence, if $Y$ is nonempty, then it contains a site $y$ with at least $d$ infected neighbors (consider the rightmost vertices in direction 1, then the rightmost of those in direction 2, and so on). But then $y \in [A]$, so $y \notin Y$. This contradiction implies that $Y$ is empty, as claimed.

Finally, note that $\mathbb{P}(B \text{ is inactive}) \leq \delta$, and so we may bound the probability of the existence of a component of inactive blocks which spans two opposite sides of $[n]^d$ using standard techniques from percolation theory. □

Lemma 11 has the following important consequence.

Lemma 12. *For each $2 \leq r \leq d \in \mathbb{N}$ and $\ell \in \mathbb{N}_0$, there exists a constant $\delta' = \delta'(d, \ell, r) > 0$ such that the following holds. Let $n, m \in \mathbb{N}$, let $\varepsilon, p > 0$ and let $A \in \mathrm{Bin}(C^*([n+1]^d \times [2]^\ell, r), p)$. Suppose that*

$$P(m, d-i, \ell+i, r-i, p) \geq 1 - \delta'$$

*for each $1 \leq i \leq r-2$, and that $M/m$ is sufficiently large. Then,*

$$\mathbb{P}([n]^d \times \mathbf{1}^\ell \subset [A \cup ([M]^d \times \mathbf{1}^\ell)]) \geq 1 - \varepsilon,$$

*and so, in particular,*

$$P(n, d, \ell, r, p) \geq (1 - \varepsilon)P(M, d, \ell, r, p)$$

*whenever $M \leq n \in \mathbb{N}$.*

Proof. For each $t \in \mathbb{N}$ and $S \subset [d]$, consider the set

$$M_t(S) = \{\mathbf{x} \in \mathbb{N}^d : x_i \in \{t+1, t+2\} \text{ if } i \in S \text{ and } x_i \in [t] \text{ if } i \notin S\},$$

and let $M_t^*(S) = \{\mathbf{x} \in M_t(S) : x_i = t+1 \text{ if } i \in S\}$. For each $S \subset [d]$, define a bootstrap structure $C_t(S)$ on $M_t(S) \times [2]^\ell$ by giving threshold $\max\{r - |S|, 0\}$ to the elements of $M_t^*(S) \times \mathbf{1}^\ell$ and threshold $r + \ell$ to the others. Observe that, when $r - |S| \geq 2$, this structure is a copy of $C^*([t]^{d-|S|} \times [2]^{\ell+|S|}, r - |S|)$.

Claim. *Suppose that $[m]^d \times \mathbf{1}^\ell \subset [A]$, and that $C_t(S)$ is internally semi-spanned for each $\varnothing \neq S \subset [d]$ and each $m \leq t \leq n-1$. Then $[n]^d \times \mathbf{1}^\ell \subset [A]$.*

Proof. It is sufficient to prove the claim for $n = m+1$. First let $j \in [d]$, and consider an element $\mathbf{x} \in M_m^*(\{j\}) \times \mathbf{1}^\ell$. It has a neighbor in $[m]^d \times \mathbf{1}^\ell$. Thus, if $[m]^d \times \mathbf{1}^\ell \subset [A]$ and $C_m(\{j\})$ is internally semi-spanned, then $M_m^*(\{j\}) \times \mathbf{1}^\ell \subset [A]$.

In general, let $\partial S = \{S \setminus \{i\} : i \in [d]\}$ denote the shadow of $S$, and observe that the sets $\{M_m^*(U) : U \in \partial S\}$ are pairwise disjoint. Thus, $\mathbf{x} \in M_m^*(S) \times \mathbf{1}^\ell$ has $|S|$ neighbors in the set $\bigcup_{U \in \partial S} M_m^*(U) \times \mathbf{1}^\ell$. Hence, if $\bigcup_{U \in \partial S} M_m^*(U) \times \mathbf{1}^\ell \subset [A]$ and $C_m(S)$ is internally semi-spanned, then $M_m^*(S) \times \mathbf{1}^\ell \subset [A]$. Hence, the sets $M_m^*(S) \times \mathbf{1}^\ell$ are infected in order of increasing $|S|$. Finally, note that $[m+1]^d \times \mathbf{1}^\ell = \bigcup_S M_m^*(S) \times \mathbf{1}^\ell$.   $\square$

Now, choose $\delta(d', \ell', r')$ and $C(d', \ell', r')$ according to Lemma 11 for each $2 \leq r' \leq d' \in \mathbb{N}$ and $\ell' \in \mathbb{N}_0$, and let $\delta' = \min\{\delta(d-i, \ell+i, r-i) : i \in [r-2]\}$ and $C = \max\{C(d-i, \ell+i, r-i) : i \in [r-2]\}$. Then, using the FKG inequality and Lemma 11,

$$P(n, d, \ell, r, p)$$

$$\geq P(M, d, \ell, r, p) \prod_{t=M}^\infty \left(1 - \sum_{S \subset [d]} [1 - P(t, d-|S|, \ell+|S|, r-|S|)]\right)$$

$$\geq P(M, d, \ell, r, p) \prod_{t=M}^\infty (1 - 2^d C e^{-t/m})$$

$$\geq (1-\varepsilon) P(M, d, \ell, r, p)$$

if $M/m$ is sufficiently large, as required.   $\square$

We can now prove Theorem 5.

PROOF OF THEOREM 5.   The proof is by induction on $r$. The theorem holds for $r = 2$ by Lemma 9, so let $r \geq 3$ and assume the result holds for all smaller values of $r$, for all values of $d \geq r$, $\ell \geq 0$ and $\varepsilon > 0$. We shall fill the set $[n]^d$ in three steps: First we use Lemma 10 to fill a cube of sidelength $M \approx (\log n)^{1-\varepsilon}$; then we use Lemma 12 to fill a cube of sidelength $N = (\log n)^3$; finally we show that such an internally spanned cube exists somewhere in $[n]^d$ with high probability, and that this cube grows to fill all of $[n]^d$.

Let $d$, $\ell$, $r$ and $\varepsilon$ be as described, let $n$ be sufficiently large, and let

$$p \geq \left(\frac{\lambda(d+\ell, \ell+r) + \varepsilon}{\log_{(r-1)} n}\right)^{d-r+1}.$$

Let $\delta = \delta(d, \ell, r, \varepsilon) > 0$ be sufficiently small, and let $m$ be defined by

$$\log_{(r-2)} m = (1 - 2\delta) \log_{(r-1)} n,$$

let $M$ be defined by $\log_{(r-2)} M = (1-\delta) \log_{(r-1)} n$, and let $N = (\log n)^3$. Note that $M/m \to \infty$ as $n \to \infty$.

First we give a lower bound on the probability that $[M]^d \times [2]^\ell$ is semi-spanned.

CLAIM 1.   $P(M, d, \ell, r, p) \gg \frac{1}{n}$ *as* $n \to \infty$.

PROOF.   When $r = 3$ this follows from Lemma 10. Indeed, note that $\log M = (1 - \delta) \log \log n$, so $p^{-2d} \le (\log \log n)^{2d^2} \le M^\delta$, and thus,

$$P(M, d, \ell, 3, p) \ge \exp(-M^{1+\delta}) = \exp(-(\log n)^{1-\delta^2}) \gg \frac{1}{n}.$$

When $r \ge 4$ the claim is even easier, since $d \log M \le d(\log \log n)^{1-\delta} \ll (1-\delta) \log \log n$, so

$$P(M, d, \ell, r, p) \ge p^{M^d} \gg \exp(-\log(1/p)(\log n)^{1-\delta}) \gg \frac{1}{n}. \qquad \square$$

Next we apply the induction hypothesis to show that semi-percolation is likely to occur on the sides of $[m]^d \times [2]^\ell$.

CLAIM 2.   $P(m, d-i, \ell+i, r-i, p) \to 1$ *as* $n \to \infty$ *for every* $1 \le i \le r-2$.

PROOF.   Observe that

$$p \ge \left( \frac{(1-2\delta)(\lambda(d+\ell, \ell+r) + \varepsilon)}{\log_{(r-2)} m} \right)^{d-r+1} \ge \left( \frac{\lambda(d+\ell, \ell+r) + \delta}{\log_{(r-i-1)} m} \right)^{d-r+1}$$

for each $i \in [r-2]$, if $\delta$ is sufficiently small. Thus, for each $i \in [r-2]$,

$$P(m, d-i, \ell+i, r-i, p) \to 1$$

as $n \to \infty$, by the induction hypothesis.   $\square$

Claim 2 allows us to apply Lemma 12. Combining this with Claim 1, we obtain

$$P(N, d, \ell, r, p) \ge \frac{1}{2} P(M, d, \ell, r, p) \ge \frac{1}{n}$$

if $n$ is sufficiently large. It follows, as there are $(n/N)^d \gg n$ pairwise disjoint cubes, that, with high probability, there exists a cuboid $C \times \mathbf{1}^\ell \subset [A]$ of size $[N]^d \times \mathbf{1}^\ell$ somewhere in $[n]^d \times [2]^\ell$. Applying Lemma 12 and the FKG inequality once more, we obtain

$$P(n-1, d, \ell, r, p) \ge (1 - o(1)) \mathbb{P}([n-1]^d \times \mathbf{1}^\ell \subset [A \cup (C \times \mathbf{1}^\ell)]) \to 1$$

as $n \to \infty$. One final application of the induction hypothesis to the sets $\{x \in [n]^d \times [2]^\ell : x_j = n\}$ now gives $\mathbb{P}([n]^d \times \mathbf{1}^\ell \subset [A]) \to 1$ as $n \to \infty$, as required. $\square$

**4. Percolation on $C([n]^2 \times [k], 2)$.**  In this section we shall prove the following theorem and corollary, from which Theorem 2 follows.

THEOREM 13.  *For every $\varepsilon > 0$, there exists $B_0 > 0$ and $k_0 : \mathbb{R}^+ \to \mathbb{R}^+$ such that the following holds for all $B \geq B_0$ and $k \geq k_0(B) \geq 3$. Let $p > 0$ be sufficiently small, let $R \subset C(B/p, k)$ be a rectangle with $\mathrm{long}(R) = B/p$, and let the elements of $A \subset R$ be chosen independently at random with probability $p$. Then*

$$\mathbb{P}(R \in \langle A \rangle) \leq \exp\left(-\frac{2\lambda(3,3) - \varepsilon}{p}\right).$$

The following corollary is the technical statement which we shall need in Section 5.

COROLLARY 14.  *For every $\varepsilon > 0$, there exist $B, k_0 > 0$ such that if $k \geq k_0$, $n$ is sufficiently large, and the elements of $A \subset C(n, k)$ are chosen independently at random with probability $p = \frac{\lambda(3,3) - \varepsilon}{\log n}$, then*

$$\mathbb{P}(\mathrm{long}(R) \geq B \log n \ for \ some \ R \in \langle A \rangle) \leq n^{-\varepsilon}.$$

Our proof will be similar in structure to that given by Holroyd [19] in the 2-dimensional case; however, the proof does not follow from that of [19] in a straightforward way. In fact, even our notion of a "hierarchy" is different, and this makes "crossing a rectangle" somewhat harder. In Sections 4.1 and 4.2 we make the necessary definitions and deal with the resulting technical problems. Finally, in Section 4.3 we sketch how the method of [19] may be used to complete the proof.

One of the important ideas of Holroyd was that the bootstrap process in a $(B/p) \times (B/p)$ rectangle may be broken up into a *bounded* number of steps, each step being either the appearance of a small internally filled rectangle (a "seed"), the growth of a rectangle sideways by $\varepsilon/p$, or the combination of two (not too small) rectangles into a larger one. Moreover, and crucially, these steps are caused by *disjoint* sets of active sites, so, having bounded the probability of each step, the probability of a particular "hierarchy" of rectangles may be bounded from above using the van den Berg–Kesten Lemma. The point is that there are either many "sideways steps" or many seeds.

In our case the situation is a little more complicated, and we therefore have to define the hierarchy slightly differently (see Section 4.1), using the concept of internal spanning defined in Section 2. It is then somewhat trickier to bound the probability that a rectangle grows sideways: we do this is Section 4.2. Bounding the probability of a seed appearing is easy, as in [19].

We remark here, for ease of reference, that there will be various constants which appear in the proof, which will depend on each other, but *not* on $p$.

These will be chosen in the order first $B$ (for "big"), then $\delta$, $k$ and $Z$ (for "seed") together, and finally $T$ (for "tiny"), and will satisfy

$$T \ll \delta, Z \ll 1 \ll B \ll k.$$

In particular, we shall need that $\delta \leq \delta(B)$, $k \geq k(B, \delta)$ and $T \leq T(Z, k, \delta)$ in Lemmas 21 and 28, that $\delta \leq \delta(Z)$ in Lemma 28, and that $Z \leq Z(B, k)$ in Lemma 31. Fortunately all of these inequalities can be satisfied simultaneously, as we shall see.

4.1. *Hierarchies.* The purpose of this subsection is to prove Lemma 20, below, which gives us our fundamental bound on the probability that $A$ percolates. (Note that we are now referring to full percolation, as opposed to the semi-percolation studied in the previous section.) In order to state the lemma, we shall need to define what we mean by a "good and satisfied hierarchy" of a rectangle $R$. In this subsection we shall work in $C([n]^d \times [k]^\ell, r)$, since the proofs carry over to the general case in a very natural way. We shall assume throughout that $k \geq 3$, although in fact our proofs also work in the case $k = 2$.

We begin by describing the algorithm by which we infect the sites of $C([n]^d \times [k]^\ell, r)$. It is slightly more complicated than the algorithm used in [19], and may seem slightly unnatural at first. Defining a hierarchy in this way seems to be necessary, however, and is perhaps the most crucial new idea in this paper.

Main Algorithm. At each step of the algorithm we have a collection of rectangles $R_1, \ldots, R_m$ and a collection of *disjoint* sets $A_1, \ldots, A_m \subset A$ such that $\langle A_i \rangle = R_i$ for each $i \in [m]$, that is, $\Pi([A_i])$ is connected and $R_i$ is the smallest rectangle containing $[A_i]$.

We begin by letting $m = |A|$, and partitioning $A$ into single elements, that is, $|A_i| = 1$ for each $i \in [m]$. Thus, each $R_i = \langle A_i \rangle$ is a $1 \times \cdots \times 1$ rectangle (considered in $[n]^d$). At each step, we perform one of the following operations:

(a) If $\Pi([A_i]) \cup \Pi([A_j])$ is connected, then we replace

$$(R_i, A_i) \text{ and } (R_j, A_j) \text{ by } (\langle A_i \cup A_j \rangle, A_i \cup A_j).$$

(b) If $2 \leq t \leq r + \ell$ is minimal such that $[A_{j(1)} \cup \cdots \cup A_{j(t)}] \neq [A_{j(1)}] \cup \cdots \cup [A_{j(t)}]$, then we replace the collection

$$\{(R_{j(1)}, A_{j(1)}), \ldots, (R_{j(t)}, A_{j(t)})\}$$

$$\text{by } (\langle A_{j(1)} \cup \cdots \cup A_{j(t)} \rangle, A_{j(1)} \cup \cdots \cup A_{j(t)}).$$

When neither of the operations is possible, or when $m = 1$, we stop and output the collection $\{R_1, \ldots, R_m\}$.

We claim that, after each step, the sets $\{A_j : j \in [m]\}$ are pairwise disjoint, and $\Pi([A_j])$ is connected for each $j \in [m]$, as required. For (a) this is obvious; for (b) it follows because $t$ was chosen to be minimal, and so one of the elements of $[A_{j(1)} \cup \cdots \cup A_{j(t)}] \setminus ([A_{j(1)}] \cup \cdots \cup [A_{j(t)}])$ connects the components of $\Pi([A_{j(1)}] \cup \cdots \cup [A_{j(t)}])$.

We make the following observation about the algorithm above.

OBSERVATION 15.   For any $A \subset C([n]^d \times [k]^\ell, r)$, the output of the Main Algorithm is $\langle A \rangle$.

Moreover, the diameter of the largest rectangle at most doubles at each step. Thus, we have the following two key lemmas from [19]. (We remark that the first was originally proved for $[n]^d$ in [2], and a version of the second for the hypercube was independently proved in [3].)

LEMMA 16.   Let $A \subset C([n]^d \times [k]^\ell, r)$. If $1 \le L \le \mathrm{diam}([A])$, then there exists a rectangle $R$, internally spanned by $A$, with

$$L \le \mathrm{long}(R) \le 2L.$$

PROOF.   Run the Main Algorithm for $A$. At some point along the way the required rectangle must have been created.   □

LEMMA 17.   Let $R \subset C([n]^d \times [k]^\ell, r)$ be a rectangle, and suppose that $R \in \langle A \cap R \rangle$. Then, for some $2 \le t \le r + \ell$, there exist disjoint nonempty sets $A_1, \ldots, A_t \subset A$, and rectangles $U_1, \ldots, U_t$, such that $\langle A_i \rangle = U_i \ne R$ for each $i \in [t]$, and $\langle A_1 \cup \cdots \cup A_t \rangle = R$.

REMARK 3.   Note that we prove that $\langle A_1 \cup \cdots \cup A_t \rangle = R$, not just that $\langle U_1 \cup \cdots \cup U_t \rangle = R$ (as in previous versions of the lemma, see [3, 19]). This subtlety will be important in the proof of Lemma 18 below.

PROOF OF LEMMA 17.   Let $A' \subset A \cap R$ be minimal such that $\langle A' \rangle = R$; such a set must exist since $R \in \langle A \cap R \rangle$. Run the Main Algorithm for $A'$ up until the penultimate step. Whether the last step is of Type (a) or Type (b) we obtain, for some $2 \le t \le r + \ell$, disjoint nonempty sets $A_1, \ldots, A_t$, as required. Indeed, $\langle A_i \rangle$ is a rectangle for each $i \in [t]$ by the definition of the algorithm, and $\langle A_i \rangle \ne R$ since $A'$ was chosen to be minimal. Finally, $\langle A_1 \cup \cdots \cup A_t \rangle = R$ by Observation 15 since $\langle A' \rangle = R$.   □

We need one more important definition.

DEFINITION. Given two rectangles $R \subset R'$, let $D(R, R')$ denote the event that

$$R' \in \langle (A \cup R) \cap R' \rangle,$$

that is, the event that $R'$ is internally spanned by $A \cup R$.

Note that the event $D(R, R')$ depends only on the set $A \cap (R' \setminus R)$. Let

$$P_p(R, R') := \mathbb{P}(D(R, R')|A \in \mathrm{Bin}(R', p)).$$

DEFINITION. Let $R$ be a rectangle in $C([n]^d \times [k]^\ell, r)$, and let $p > 0$. A *hierarchy* $\mathcal{H}$ of $R$ is an oriented rooted tree $G_\mathcal{H}$, with all edges oriented away from the root ("downward"), together with a collection of rectangles $\{R_u : u \in V(G_\mathcal{H})\}$, $R_u \subset C([n]^d \times [k]^\ell, r)$, one for each vertex of $G_\mathcal{H}$, satisfying the following criteria:

(a) The root of $G_\mathcal{H}$ corresponds to $R$.
(b) Each vertex has at most $r + \ell$ neighbors below it.
(c) If $u \to v$ in $G_\mathcal{H}$, then $R_u \supset R_v$.
(d) If $\vec{\Gamma}(u) = \{v_1, \ldots, v_t\}$ and $t \geq 2$, then $\langle R_{v_1} \cup \cdots \cup R_{v_t} \rangle = R_u$.

A hierarchy is *good* for $(T, Z, p) \in \mathbb{R}^3$ if:

(e) If $\vec{\Gamma}(u) = \{v\}$ and $|\vec{\Gamma}(v)| = 1$, then $\phi(R_u) - \phi(R_v) \in [T/p, 2T/p]$.
(f) If $\vec{\Gamma}(u) = \{v\}$ and $|\vec{\Gamma}(v)| \neq 1$, then $\phi(R_u) - \phi(R_v) \leq 2T/p$.
(g) If $|\vec{\Gamma}(u)| \geq 2$ and $v \in \vec{\Gamma}(u)$, then $\phi(R_u) - \phi(R_v) \geq T/p$.
(h) If $u$ is a leaf, then $\mathrm{short}(R_u) \leq Z/p$.
(i) If $u$ is not a leaf, then $\mathrm{short}(R_u) > Z/p$.

A hierarchy is *satisfied* by $A$ if the following events all occur *disjointly*:

(j) $R_u$ is internally spanned by $A$ whenever $R_u$ is a seed (i.e., $u$ is a leaf).
(k) $D(R_v, R_u)$ whenever $\vec{\Gamma}(u) = \{v\}$.

The next lemma tells us that every internally spanned rectangle $R$ has a good and satisfied hierarchy. See also Proposition 32 of [19].

LEMMA 18. *Let* $A \subset C([n]^d \times [k]^\ell, r)$, *let* $T, Z > p > 0$, *and let* $R \subset C([n]^d \times [k]^\ell, r)$ *be a rectangle. Suppose that* $A$ *internally spans* $R$. *Then there exists a good and satisfied hierarchy of* $R$.

PROOF. The lemma follows by an easy induction on $\phi(R)$. First note that the result is immediate if $\mathrm{short}(R) \leq Z/p$, by choosing the hierarchy with one element. Thus, in particular, the result holds if $\phi(R) \leq 2Z/p$.

So let $\mathrm{short}(R) > Z/p$, and apply Lemma 17. For some $2 \le t(1) \le r + \ell$, we obtain disjoint sets $A_1^{(1)}, \ldots, A_{t(1)}^{(1)} \subset A$ and rectangles $U_1^{(1)}, \ldots, U_{t(1)}^{(1)}$, such that $\langle A_i^{(1)} \rangle = U_i^{(1)} \ne R$ for each $i \in [t(1)]$, and $\langle A_1^{(1)} \cup \cdots \cup A_{t(1)}^{(1)} \rangle = R$. Choose one of the rectangles $U_{j(1)}^{(1)}$ with

$$\phi(U_{j(1)}^{(1)}) = \max\{\phi(U_i^{(1)}) : i \in [t]\},$$

and let $S_1 = U_{j(1)}^{(1)}$ and $A_1 = A_{j(1)}^{(1)}$. Note that $\langle A_1 \rangle = S_1$, and that $\phi(S_1) < \phi(R)$.

Now apply Lemma 17 to the rectangle $S_1$ to get disjoint sets $A_1^{(2)}, \ldots, A_{t(2)}^{(2)} \subset A_1$ and rectangles $U_1^{(2)}, \ldots, U_{t(2)}^{(2)}$, and hence a pair $(S_2, A_2)$ with $\langle A_2 \rangle = S_2 = U_{j(2)}^{(2)}$ as before. Repeat until one of the following occurs for some $m \in \mathbb{N}$:

(a) $\phi(R) - \phi(S_m) \in [T/p, 2T/p]$,

(b) $\phi(R) - \phi(S_m) \ge 2T/p$,

(c) $\mathrm{short}(S_m) \le Z/p$.

Note that at least one of these must occur eventually, since $\phi(S_{t+1}) \le \phi(S_t) - 1$ for all $t \in \mathbb{N}$. There are four cases to consider:

*Case* 1: $\phi(R) - \phi(S_m) \in [T/p, 2T/p]$. By induction, there exists a good [for $(T, Z, p)$] and satisfied (by $A_m$) hierarchy $\mathcal{H}'$ of $S_m$. We create a good and satisfied hierarchy $\mathcal{H}$ of $R$ by adding a new root vertex, with a single neighbor (the root vertex of $\mathcal{H}'$). It is easy to see that $\mathcal{H}$ is a good hierarchy for $(T, Z, p)$; it is satisfied by $A$ because the set $A \setminus A_m$ is a witness set for the event $D(S_m, R)$, since $A_m \subset S_m$ and $R \in \langle A \rangle$.

*Case* 2: $\phi(R) - \phi(S_1) \ge 2T/p$. There exist good and satisfied hierarchies $\mathcal{H}_1, \ldots, \mathcal{H}_{t(1)}$ for $U_1^{(1)}, \ldots, U_{t(1)}^{(1)}$ respectively, where $\mathcal{H}_i$ is satisfied by $A_i^{(1)}$ for each $i \in [t(1)]$. We obtain a good and satisfied hierarchy $\mathcal{H}$ for $R$ by adding a new root vertex, with $t(1)$ neighbors [the root vertices of $\mathcal{H}_1, \ldots, \mathcal{H}_{t(1)}$]. This hierarchy is clearly satisfied by $A$; it is good because $\phi(R) - \phi(U_i^{(1)}) \ge \phi(R) - \phi(S_1) \ge 2T/p$ for each $i \in [t(1)]$.

*Case* 3: $\phi(R) - \phi(S_m) \ge 2T/p$ for some $m \ge 2$. Since this is the first $m$ for which one of (a), (b) and (c) holds, it follows that $\phi(R) - \phi(S_{m-1}) < T/p$ and $\mathrm{short}(S_{m-1}) > Z/p$. Note that therefore $\phi(S_{m-1}) - \phi(U_i^{(m)}) \ge \phi(S_{m-1}) - \phi(S_m) \ge T/p$ for each $i \in [t(m)]$.

Let $\mathcal{H}_1, \ldots, \mathcal{H}_{t(m)}$ be good and satisfied hierarchies for $U_1^{(m)}, \ldots, U_{t(m)}^{(m)}$, respectively. Define $\mathcal{H}$ by adding two new vertices: a new root vertex $u$, with one neighbor $v$, which in turn has $t(m)$ other neighbors [the root vertices of $\mathcal{H}_1, \ldots, \mathcal{H}_{t(m)}$]. Let $R_v = S_{m-1}$, and observe that $\mathcal{H}$ is good and satisfied.

*Case* 4: short$(S_m) \leq Z/p$, but $\phi(R) - \phi(S_m) < T/p$. Let $\mathcal{H}_1$ be a good and satisfied hierarchy for $S_m$ (i.e., a single vertex), and form $\mathcal{H}$ by adding a root vertex to $\mathcal{H}_1$. It is easy to see that $\mathcal{H}$ is a good and satisfied hierarchy for $R$. $\quad\square$

Given $T, Z > p > 0$, let $\mathcal{H}(R, T, Z, p)$ denote the collection of hierarchies for $R$ which are good for the triple $(T, Z, p)$. The next lemma makes the crucial observation that there are only "few" possible hierarchies.

LEMMA 19. *Let* $B, p > 0$, *let* $R \subset C([n]^d \times [k]^\ell, r)$ *with* long $(R) \leq B/p$, *and let* $T, Z > p > 0$. *Then there exists a constant* $M = M(B, T, d, \ell, r)$ *such that*
$$|\mathcal{H}(R, T, Z, p)| \leq M p^{-M}.$$

PROOF. Let $\mathcal{H}$ be a hierarchy in $\mathcal{H}(R, T, Z, p)$, and consider a path from the root of $\mathcal{H}$ to a leaf. By properties (e) and (g), out of every two consecutive steps (not including the last), there is one which corresponds to a decrease in $\phi(R)$ of at least $T/p$. These add up to at most $dB/p$, and thus, the tree $G$ of $\mathcal{H}$ has depth at most $2dB/T + 1$. It also has maximal out-degree at most $r + \ell$ by property (b). Each such tree has at most $V = 2(r + \ell)^{2dB/T+1}$ vertices, and there are thus at most $(r + \ell + 1)^V$ such trees.

Now, each rectangle may be chosen in at most $(B/p)^{2d}$ ways, and so there are at most $(B/p)^{2dV}$ ways of choosing the rectangles. Thus, there are at most
$$(r + \ell + 1)^V (B/p)^{2dV} \leq (B(r + \ell + 1))^{2dV} p^{-2dV} \leq M p^{-M}$$
possible hierarchies of $R$, where $M = (B(r + \ell + 1))^{2dV}$. $\quad\square$

We are ready to prove the main lemma of this section. It gives us our basic bound on the probability that $A$ internally spans $R$. Recall that $P_p(R)$ denotes the probability that a rectangle $R$ is spanned by a set $A \in \mathrm{Bin}(R, p)$.

LEMMA 20. *Let* $R$ *be a rectangle in* $C([n]^d \times [k]^\ell, r)$, $T, Z > p > 0$, *and* $A \in \mathrm{Bin}(R, p)$. *Then*
$$\mathbb{P}(R \in \langle A \rangle) \leq \sum_{\mathcal{H} \in \mathcal{H}(R, T, Z, p)} \left( \prod_{\vec{\Gamma}(u) = \{v\}} P_p(R_v, R_u) \right) \prod_{\text{seeds } u} P_p(R_u).$$

PROOF. Suppose $R \in \langle A \rangle$. Then, by Lemma 18, there exists a good and satisfied hierarchy $\mathcal{H}$ for $V$. But, by the van den Berg–Kesten Lemma, the probability that $\mathcal{H}$ is satisfied by $A$ is at most the product of the probability of the following events:

(j) $R_u$ is internally spanned by $A$ whenever $R_u$ is a seed (i.e., $u$ is a leaf), and

(k) $D(R_v, R_u)$ whenever $\vec{\Gamma}(u) = \{v\}$

since these occur disjointly. The result now follows by taking the union bound over all possible hierarchies. $\square$

4.2. *Crossing a rectangle.* We now return to the three-dimensional case, $C(n, k)$. We begin by defining what we mean by crossing a rectangle in $C(n, k)$. Our definition is a generalization of that for $[n]^2$ in [19].

A path from left to right across a rectangle $R = [(a, b), (c, d)] \subset C(n, k)$ is a path from a point in the set $\{(x, y, z) \in R : x = a\}$ to a point in the set $\{(x, y, z) \in R : x = c\}$.

DEFINITION.   A rectangle $R = [(a, b), (c, d)] \subset C(n, k)$ is said to be *left-to-right crossed* (or just crossed) by $A \subset C(n, k)$ if the set $A \cap R$ has the following property: let

$$A' := (A \cap R) \cup \{(x, y, z) : x \leq a - 1\}.$$

Then there is path in $[A']$ from left to right across $R$.

We write $H^{\to}(R)$ for this event, and define $H^{\leftarrow}(R)$ (right-to-left), $H^{\downarrow}(R)$ (top-to-bottom) and $H^{\uparrow}(R)$ (bottom-to-top) crossing of $R$ similarly. (Here "top-to-bottom," e.g., means from a larger to a smaller second coordinate.) In $[n]^2$ crossing a rectangle is simple; one simply has to avoid "double gaps." In $C(n, k)$ more things can go wrong, so we begin by bounding the event $H^{\to}(R)$ (the others follow by symmetry). In fact, and with foresight, we shall bound from above the function

$$h(R, \ell) := \max_{W \subset R, |W| \leq \ell} \{\mathbb{P}_p(R \text{ is left-to-right crossed by } A | W \subset A)\},$$

where we write $\mathbb{P}_p$ to mean $A \in \text{Bin}(R, p)$. Note that $\mathbb{P}_p(H^{\to}(R)) = h(R, 0)$.

Recall the definition (1) of $\beta(u) := \beta_2(u)$ from the Introduction. In particular, note that it satisfies

$$\beta(u)^2 = (2u - u^2)\beta(u) + u(1 - u)^2.$$

The following lemma is the key (new) step in the proof of Theorem 13.

LEMMA 21.   *Let $B > 0$, $\delta_0(B) > \delta > 0$, and $k \geq 6e^{6B} \log(1/\delta) + 2$. There exists a constant $T = T(k, \delta) > 0$ such that the following holds. Let $p > 0$ be sufficiently small, and let $R$ be a rectangle in $C(n, k)$, with $\dim(R) = (s, m)$, where $m \leq B/p$ and $s \leq T/p$. Then, for any $\ell \in \mathbb{N}$ with $2\ell \leq s$,*

$$\beta(u)^{s+1} \leq \mathbb{P}_p(H^{\to}(R)) \leq h(R, \ell) \leq m^3(\beta(u) + \delta)^{s - 4\ell},$$

*where $u = 1 - (1 - p)^m = 1 - e^{-qm}$.*

REMARK 4. To understand this lemma, the reader should think of $B$ and $k$ as large, and of $\beta(u)$, $\delta$ and $T$ as constants, with $\delta$ smaller than $1 - \beta(u) \in (0,1)$, and $T$ *much* smaller than $\delta$. As we shall see later, the error terms $m^3$, $\delta$ and $4\ell$ on the right-hand side do not matter much, and so the lemma gives an essentially sharp upper bound on $h(R, \ell)$.

For the sake of simplicity, we shall assume that $k$ is even in the proof of Lemma 21; the proof for $k$ odd is the same. Thus, we now replace $k$ by $2k$, and throughout the remainder of this section we let $R = [(1,1),(s,m)] \subset C(n, 2k)$ be a rectangle as in Lemma 21, and assume that the set $\{(0, y, z) : (1, y, z) \in R\} \subset A$. We begin by defining some events which depend on the set $A \cap R \in \text{Bin}(R, p)$, which we shall call *blockers*, *savers* and *last chances* (see Figure 1).

In Figure 1 the top left point of $R$ is $(1,1,1)$, the $x$-axis runs left-to-right, the $z$-axis top to bottom and the $y$-axis into the page. Thus, the top and bottom surfaces of $R$ (in the figure) have threshold 2, and every vertex not on one of these surfaces has threshold 3. The shaded area on the left denotes the set $A' \setminus R$ of "previously infected" sites.

We begin by defining some sets, which are each $1 \times m \times 1$ columns (going into the page). For each $i \in [s]$ and $j \in [k]$, let

$$M_i(j) := \{(x, y, z) \in R : x = i \text{ and } z = j\},$$

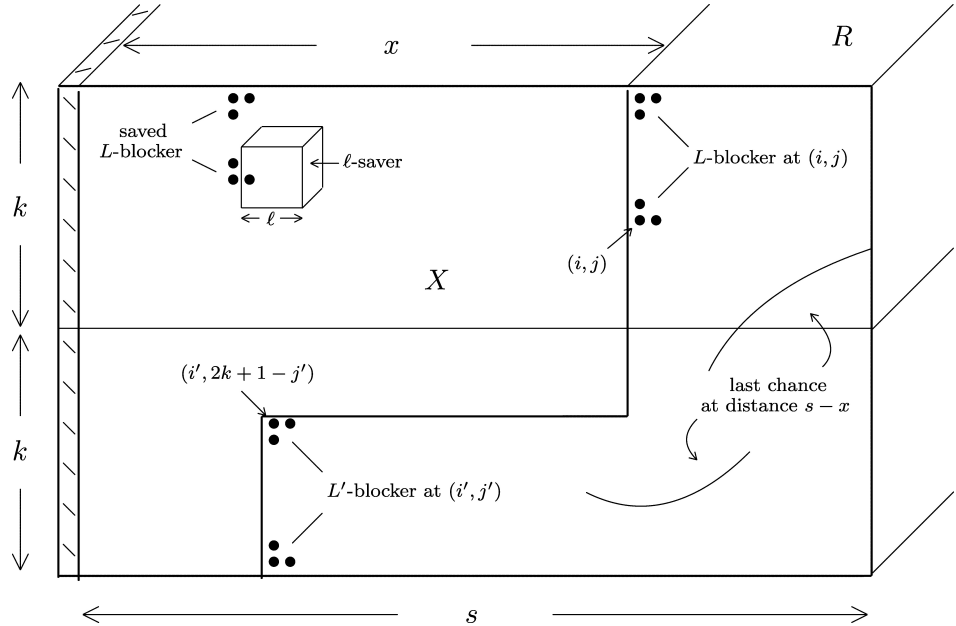

FIG. 1. *Crossing $R$.*

and let $M_i'(j) := M_i(2k + 1 - j)$.

DEFINITION. An *L-blocker* occurs at point $(i, j)$ if the set

$$M_i(1) \cup M_i(2) \cup M_{i+1}(1) \cup M_i(j) \cup M_i(j - 1) \cup M_{i+1}(j)$$

is empty (i.e., contains no element of $A$), where $i \in [s - 1]$ and $2 \leq j \leq k$. Define an *L'-blocker* similarly, for the sets $M_i'(j)$.

The $L$-blockers act like the $L$-gaps of Section 3 (with $\ell = 1$). The idea is that, if $k$ is sufficiently large, the probability that there is an $L$-blocker at point $(i, j)$ for *some* $j \in [k]$ is about the same as the probability of an $L$-gap.

For each $1 \leq i \leq j \leq s$, let $R[i, j] = \{(a, b, c) \in R : a \in [i, j]\}$, and let $R[i] = R[i, i]$. Moreover, let $R^+ = \{(a, b, c) \in R : c \leq k\}$ and $R^- = R \setminus R^+$, and define $R^+[i, j]$ and $R^-[i, j]$ accordingly.

A *double gap* in a cuboid $C = [a_1] \times [a_2] \times [a_3]$ is a pair of empty adjacent planes in $C$, that is, a pair $(i, j)$, with $0 \leq j \leq a_i$, such that $\{(x_1, x_2, x_3) \in A \cap C : x_i \in \{j, j + 1\}\} = \varnothing$. Note that this definition includes the case where just the face of the cuboid is empty.

Our next definition deals with the possibility that, although $M_i(j)$ may be empty (in $A$), it may contain some element of $[A]$.

DEFINITION. An *$\ell$-saver* of $M_i(j)$ is a cuboid $C \subset R^+[i, i + \ell - 1]$ whose left face intersects $M_i(j)$, whose right face intersects $M_{i+\ell-1}(j)$, and which has no double gap.

An *L*-blocker at point $(i, j)$ is said to be *saved* (or *$\ell$-saved*) if there exists an $\ell$-saver of one of $M_i(1)$, $M_i(2)$, $M_{i+1}(1)$, $M_i(j)$, $M_i(j - 1)$ and $M_{i+1}(j)$, for some $2 \leq \ell \leq s$. Otherwise, the $L$-blocker is *unsaved*.

We define $\ell$-savers and *L'*-blockers in $R^-$ similarly, using the sets $M_i'(j)$.

We shall show (see Lemma 26 below) that an $L$-blocker is very unlikely to be saved, and thus that the probability that there is an *unsaved* $L$-blocker at point $(i, j)$ for some $j \in [k]$ is also about the same as the probability of an $L$-gap.

The following algorithm describes a method of trying to cross $R^+$, and defines the variable $CA(R^+) \in [0, s]$.

CROSSING ALGORITHM. Set $CA(R^+) := 0$ and $x := 0$, and repeat the following steps until either $CA(R^+) \geq s$, or STOP:

1. If the set $M_{x+1}(1) \cup M_{x+1}(2)$ is occupied, then set $x := x + 1$, and go to Step 5.
2. If the set $M_{x+2}(1)$ is occupied, then set $x := x + 2$, and go to Step 5.

3. If $M_{x+1}(1) \cup M_{x+1}(2) \cup M_{x+2}(1)$ is empty, but there is no $L$-blocker at point $(x+1, j)$ for any $2 \le j \le k$, then set $x := x+2$, and go to Step 5.

4. Otherwise, let $u \in [k]$ be minimal such that there is an $L$-blocker at point $(x+1, u)$. Set $\ell := 2$, and repeat the following steps until $x + \ell > s$:

   (a) If there exists an $\ell$-saver of the $L$-blocker at point $(x+1, u)$, then set $x := x + \ell + 1$ and go to Step 5.

   (b) Otherwise,

      (i) If $x + \ell \ge s$, then STOP.

      (ii) Set $\ell := \ell + 1$ and go back to Step 4(a).

5. Set $CA(R^+) := x$ and go back to Step 1.

REMARK 5.    Note that if an $L$-blocker is saved by an $\ell$-saver, we "give away for free" the next $\ell + 1$ columns, about which we now have "positive" information (that some cuboid has no double gaps), in order to preserve independence. However, we pay a price for this: the Crossing Algorithm is not monotone. [For example, adding an infected site in $M_i(1)$ is unhelpful if we would otherwise use a 10-saver of $M_i(5)$.] It is for this reason that we will need to use Reimer's Theorem.

Using the Crossing Algorithm, we come to the definition we shall use.

DEFINITION.    Say that $R^+$ is $L$-*crossed* up to the point $x \in [s]$ if $CA(R) \ge x$. Define $L'$-crossing of $R^-$ similarly, using $L'$-blockers and the sets $M_i'(j)$ in the Crossing Algorithm.

If either $R^+$ is $L$-crossed or $R^-$ is $L'$-crossed up to $x$, then say that $R$ is *unblocked* up to $x$.

We shall use the following properties of $L$-crossing.

LEMMA 22.    *For any* $x, y \ge 0$,

$$\mathbb{P}_p(R \text{ is } L\text{-crossed up to } x+y \mid R[1,x] \subset A) \le \mathbb{P}_p(R \text{ is } L\text{-crossed up to } y).$$

PROOF.    None of the elements of $A \cap R[1,x]$ are useful in crossing from $x$ to $x+y$. The inequality comes from the fact that there is less space to the right in which to find $\ell$-savers (only $s - x$ instead of $s$).  □

LEMMA 23.    *Let* $x \in \mathbb{N}$ *be maximal such that* $R$ *is unblocked up to* $x$. *Suppose* $x \le s - 1$. *Then, for some* $(u, v)$ *and* $(u', v')$, *with* $\max\{u, u'\} = x + 1$ *and* $v, v' \in [k]$,

   (a) *there is an unsaved $L$-blocker at point* $(u, v)$, *and*

   (b) *there is an unsaved $L'$-blocker at point* $(u', v')$.

*Moreover, suppose $u' \leq u$, and let*

$$X = \{(a, b, c) \in R : a < u', \text{ or } a < u \text{ and } c \leq 2k - v'\}.$$

*Then there is a witness set in $X$ for the event "$R$ is unblocked up to $x$."*

PROOF.   This follows immediately from the Crossing Algorithm.   □

We need one more definition. Let $[A \cap R]_2$ denote the closure of the set $A \cap R$ under the 2-neighbor rule.

DEFINITION.   A *last chance* at distance $y \geq 0$ is a path in $[n]^2$, which uses only vertices from the set

$$\Pi([A \cap R]_2 \cup R[s+1]),$$

from $R[s+1]$ to $R[s-y+1]$, that is, from the boundary of $R$ on the right, to a point at distance $y$ from the boundary.

Finally our effort is rewarded: the following lemma shows why the events above are important.

LEMMA 24.   *Let $R \subset C(n, 2k)$ be a rectangle as described in Lemma 21, and let $A \subset R$. If the event $H^{\rightarrow}(R)$ occurs, then there exists some $x \leq s$ such that the events,*

(a)  *$R$ is unblocked up to $x$, and*
(b)  *there exists a last chance at distance $y = s - x$,*

*occur disjointly.*

PROOF.   Let $R = [(1, 1), (s, m)]$, and let $x \in [s]$ be maximal such that $R$ is unblocked up to $x$. If $x = s$, then we are done, since there is always a last chance at distance 0. Otherwise, by Lemma 23, there exists an unsaved $L$-blocker at $(u, v)$ and an unsaved $L'$-blocker at $(u', v')$, say, where, without loss of generality, $u' \leq u = x + 1$. Moreover, writing

$$X = \{(a, b, c) \in R : a < u', \text{ or } a < u \text{ and } c \leq 2k - v'\}$$

(see Figure 1), there is a witness set in $A \cap X$ for the event "$R$ is unblocked up to $x$." Suppose the event $H^{\rightarrow}(R)$ occurs; we claim that there is a witness set in $A \setminus X$ for the event "there exists a last chance at distance $y = s - x$."

Assume there is no such witness set, and let $A' = A \cup \{(0, b, c) : b \in [m], c \in [2k]\}$. We must show that there is no path across $R$ in $[A']$. Consider the set

$$Y = M_u(1) \cup M_u(2) \cup M_{u+1}(1) \cup M_u(v) \cup M_u(v-1) \cup M_{u+1}(v)$$

$$\cup M'_{u'}(1) \cup M'_{u'}(2) \cup M'_{u'+1}(1) \cup M'_{u'}(v') \cup M'_{u'}(v'-1) \cup M'_{u'+1}(v'),$$

where the sets $M_i(j)$ and $M'_i(j)$ are as defined above. Observe that $A \cap Y = \varnothing$, since $R$ has an $L$-blocker at point $(u, v)$ and an $L'$-blocker at point $(u', v')$.

CLAIM.  $[A'] \cap Y = \varnothing$.

PROOF.  This follows because the $L$-blocker and the $L'$-blocker are unsaved. Indeed, suppose that $[A'] \cap Y$ is nonempty, and run the bootstrap process until some element of $Y$ is infected. Let $w \in [A'] \cap Y$ be the first element of $Y$ infected by $A'$, and let $W$ denote the set of infected sites if the bootstrap process is stopped as soon as $w$ becomes infected. (To be precise, we choose an ordering $v_1, \ldots, v_t$ of $[A'] \setminus A'$ such that $|\Gamma(v_i) \cap (A' \cup \{v_1, \ldots, v_{i-1}\})| \geq r(v_i)$ for each $i \in [t]$, let $w = v_j \in Y$ with $j$ minimal, and let $W = A' \cup \{v_1, \ldots, v_j\}$.)

There are two cases to consider: either $w \in R^+$ or $w \in R^-$. Consider the sets

$$W^+ := W \cap \{(a, b, c) \in R : a \geq u \text{ and } c \leq v\}$$

and

$$W^- := W \cap \{(a, b, c) \in R : a \geq u' \text{ and } c \geq 2k + 1 - v'\}.$$

If $w \in R^+$, then let $D$ denote the connected component in $W^+$ containing $w$. If $w \in R^-$, then let $D$ denote the connected component in $W^-$ containing $w$. In both cases, let $C$ denote the smallest cuboid containing $D$.

First note that the vertices of $M_u(1) \cup M'_{u'}(1)$ have only one neighbor in $R$ outside $Y$, and the vertices of $M_u(v) \cup M'_{u'}(v')$ have only two neighbors in $R$ outside $Y$, and so $w \notin M_u(1) \cup M_u(v) \cup M'_{u'}(1) \cup M'_{u'}(v')$, since it is the first element of $Y$ infected. Next, observe that $w$ cannot lie in the right-hand edge of $C$, since $w$ would not have enough previously infected neighbors in $[A']$ to be infected itself. For example, if $w \in M_u(v - 1)$, then $w$ would have at most two previously infected neighbors, one in $X$ and one in the row above. Thus, if $C$ does not have a double gap, then it is a saver of $M_i(j)$, where $w \in M_i(j) \subset Y$.

But the blockers are unsaved, so $C$ must have a double gap, $U$. Since $C$ is the smallest cuboid containing the connected component $D$, it follows that $D$ contains some member $x \in U$. Let $x$ be the first member of $U$ to be infected; we claim that in fact $x$ must have fewer than $r(x)$ infected neighbors, a contradiction.

Indeed, $x$ has no neighbors in $D \setminus C = \varnothing$, and at most one infected neighbor in $C$, since $U$ is a double gap, and $x$ is the first member of $U$ to be infected. Moreover, $x$ has at most one neighbor in $W \setminus W^\pm$ [since $w \notin M_u(v) \cup M'_{u'}(v')$], and if $r(x) = 2$, then $x$ has no neighbors in $W \setminus W^\pm$ [since $w \notin M_u(1) \cup M'_{u'}(1)$]. It follows that $x$ has too few infected neighbors; this contradiction proves the claim.  □

Now, observe that the sites of $R \setminus (X \cup Y)$ either have threshold 3 and at most one neighbor in $X$, or have threshold 2 and no neighbors in $X$. Thus, given $[A'] \cap Y = \varnothing$, it follows that $[A']$ is contained in the set $X \cup [A \setminus X]_2$.

Recall our earlier assumption, that there is no witness set in $R \setminus X$ for a last chance at distance $y = s - x$. This means that there is no path in $\Pi([A \setminus X]_2)$ from the right-hand edge of $R$ to the set $\Pi(X)$, and so, since $[A'] \subset X \cup [A \setminus X]_2$, there is no path across the rectangle $R$ in $\Pi([A'])$. So $R$ is not crossed by $A$, contradicting our assumption that the event $H^{\rightarrow}(R)$ occurs.

We have shown that there exist witness sets for events (a) and (b) which lie in $A \cap X$ and in $A \setminus X$ respectively. Thus, the events occur disjointly, as required. $\square$

Now, define

$$a(x, \ell) = \max_{W \subset R, |W| \le \ell} \mathbb{P}_p(R^+ \text{ is } L\text{-crossed up to point } x | W \subset A).$$

Note that $a(x, \ell)$ is decreasing in $x$ and increasing in $\ell$. Most of the rest of the work of this section will be to prove the following lemma, which gives us our bound on $a(x, \ell)$.

LEMMA 25. *Let the constants $B, \delta, p > 0$ and $n, m, s, k \in \mathbb{N}$, and the rectangle $R \subset C(n, 2k)$, be as in Lemma 21, and let $u = 1 - (1-p)^m$. Then, for any $x \in [s]$ and any $\ell \in \mathbb{N}$ with $2\ell \le s$,*

$$a(x, \ell) \le (\beta(u) + \delta)^{x - 2\ell}.$$

REMARK 6. *Recall that $k \ge 6e^{6B} \log(1/\delta) + 2$, and that $s \le T/p$, where $T = T(k, \delta)$.*

In order to prove Lemma 25, we must estimate the probability that a blocker is saved. Since the proof is similar, and we shall need the result later, we shall also bound the probability that a last chance occurs. For each $y \in [s]$ each $\ell \in \mathbb{N}$, let

$$S(y, \ell) = \max_{i,j} \max_{|W| \le \ell} \{ \mathbb{P}_p(M_i(j) \text{ has a } y\text{-saver} | W \subset A) \}$$

and

$$b(y, \ell) = \max_{|W| \le \ell} \{ \mathbb{P}_p(R \text{ has a last chance at distance } y | W \subset A) \}.$$

LEMMA 26. *Let $\ell, y \in \mathbb{N}$, with $y \ge 2$, and let $n, m, s, k \in \mathbb{N}$ and the rectangle $R \subset C(n, 2k)$ be as described in Lemma 21. Then,*

(a) $S(y, \ell) \le 4m(y + 1)^3 (4kyp)^{\lceil (y+1)/2 \rceil - \ell}$,
(b) $b(y, \ell) \le 2my(6kyp)^{y/2 - \ell}$.

We shall use the following simple observation in the proof of Lemma 26.

LEMMA 27.   *Let $C \subset [n]^3$ be a cuboid with $\dim(C) = (u, v, w)$, and let $p > 0$. Let $L \subset C$ with $|L| = \ell$, and let $A \in \mathrm{Bin}(C, p)$. Then*

$$\mathbb{P}(C \text{ has no double gap} \mid L \subset A) \le (2uvp)^{\lceil (w+1)/2 \rceil - \ell}.$$

PROOF.   We partition $C$ into double slices $D_1, \ldots, D_{w'}$, where $w' = \lceil (w + 1)/2 \rceil$, by letting $D_i = \{(x, y, z) \in C : z \in \{2i-2, 2i-1\}\}$ for each $i \in [w']$. If $C$ has no double gap, then each $D_i$ is occupied, so let $W_L = \{i \in [w'] : D_i \cap L = \varnothing\}$, and let $w'' = |W_L| \ge w' - \ell$.

Let us choose, for each $i \in W_L$, an infected site $d_i \in D_i$ and let $D = \{d_i : i \in W_L\}$. We have $(2uv)^{w''}$ choices for the set $D$, and the probability that $A$ contains $D \setminus L$ is at most $p^{w''}$, as required.   $\square$

We now prove Lemma 26.

PROOF OF LEMMA 26.   Recall that $R$ is a rectangle as in Lemma 21, and let $W \subset R$ with $|W| \le \ell$. All probabilities in this proof will be conditional on the assumption that $W \subset A$. The proof in each case follows easily by counting cuboids and using Lemma 27.

Indeed, recall that a $y$-saver of $M_i(j)$ is a cuboid $C \subset R[i, i + y - 1]$ such that:

- the left face of $C$ intersects $M_i(j)$,
- the right face of $C$ intersects $M_{i+y-1}(j)$, and
- $C$ has no double gap.

Let $\mathrm{long}(C) = t \ge y$, and count cuboids. We have at most $m$ choices for the "nearmost" point in $M_i(j) \cap C$, and at most $t^3$ choices for $C$, given this point ($t$ choices in direction 2, $t^2$ choices in direction 3, and only one choice in direction 1). Note that the shorter two dimensions of $C$ are at most $y$ and $2k$ respectively. Thus, by Lemma 27, the probability that $C$ has no double gap is at most

$$(4kyp)^{\lceil (t+1)/2 \rceil - \ell}.$$

Recall that $y \le s \le T/p$, and that we may choose $T = T(k, \delta)$ as small as we like. Thus, we may assume that $4kyp$ is arbitrarily small. Hence, summing over $t$, we get

$$S(y, \ell) \le m \sum_{t=y}^m t^3 (4kyp)^{\lceil (t+1)/2 \rceil - \ell} \le 4m(y+1)^3 (4kyp)^{\lceil (y+1)/2 \rceil - \ell},$$

as claimed. In the second inequality, note that the maximum could occur at either $t = y$ or $t = y + 1$.

Next, consider $b(y, \ell)$, and recall that a last chance at distance $y$ is a path in $\Pi([A \cap R]_2 \cup R[s+1])$ from $R[s+1]$ to $R[s-y+1]$. We now have to count rectangles (in $[n]^2$), and bound the probability that each is crossed by the projection of $A$.

Indeed, suppose there is such a (shortest) path $P$ from $R[s+1]$ to $R[s-y+1]$, and consider the smallest rectangle $S \subset [n]^2$ containing the component of $\Pi([A \cap R]_2)$ which contains $P \cap R$. Then $S$ must have no double gap in $\Pi(A)$. By Lemma 27 (applied to the cuboid $S \times \{1\}$ with density $2kp$), the probability of this is at most

$$(4kup)^{\lceil (t+1)/2 \rceil - \ell},$$

where $S$ is a $u \times t$ rectangle, and $u \le t$ say. Since $P \cap R$ is a path from $R[s]$ to $R[s-y+1]$, we have $t \ge y$. Also $u \le s \le T/p$, so $4kup$ may be made arbitrarily small by an appropriate choice of $T$.

Thus, summing over all rectangles $S$, and noting that we have at most $2m$ choices for $S$ for each pair $(u,t)$, we obtain

$$b(y, \ell) \le 2m \sum_{u,t : u \le t, y \le t} (4kup)^{\lceil (t+1)/2 \rceil - \ell} \le 2my(6kyp)^{y/2 - \ell},$$

as claimed.  $\square$

Now we use Lemma 26 to prove Lemma 25, that is, to bound from above the probability that $CA(R^+) \ge x$.

Proof of Lemma 25.   Suppose that $\beta(u) + \delta < 1$ (the result is otherwise trivial). We are required to prove that, for any $x \in [s]$, any $\ell \in \mathbb{N}$ and any $W \subset R$ with $|W| \le \ell$,

$$\mathbb{P}_p(CA(R^+) \ge x | W \subset A) \le (\beta(u) + \delta)^{x - 2\ell}.$$

Note that $\mathbb{P}_p(CA(R^+) \ge x)$ depends on $s$, and in fact is increasing in $s$ (the probability of a saver existing increases with $s$). However, we shall only need the fact that $s$ is bounded from above by $T/p$, and so shall suppress this dependency on $s$.

The proof is by induction on $x + \ell$. If $x \le 2\ell$, then the result is immediate, since $a(x, \ell) \le 1$. The induction step follows easily from the following claim.

Claim.

$$\mathbb{P}_p(CA(R^+) \ge x | W \subset A)$$

$$\le \max\Big\{ a(x-2, \ell-1), (2u - u^2)a(x-1, \ell) + u(1-u)^2 a(x-2, \ell)$$

$$+ \delta^3 a(x-2, \ell) + 6 \sum_{y \ge 2} \sum_{\ell' \ge 0} S(y, \ell') a(x - y - 1, \ell - \ell') \Big\}.$$

PROOF. It follows from the Crossing Algorithm (see also Lemma 6 of Section 3) that one of the following holds:

- $M_1(1) \cup M_1(2)$ is occupied,
- $M_1(1) \cup M_1(2)$ is empty but $M_2(1)$ is occupied,
- $M_1(1) \cup M_1(2) \cup M_2(1)$ is empty but there is no $L$-blocker at $(1, j)$ for any $j \in [k]$,
- there is a $y$-saver of an $L$-blocker at point $(1, j)$ [where $j$ and $y$ are minimal, in the sense that there is no $L$-blocker at $(1, j')$ for any $j' < j$, and no $y'$-saver of the $L$-blocker at $(1, j)$ for any $y' < y$],
- $CA(R^+) = 0$.

Suppose first that $W \cap R[1, 2] \neq \varnothing$. We claim that

$$\mathbb{P}_p(CA(R^+) \geq x | W \subset A)$$
$$\leq \max \left\{ a(x - 2, \ell - 1), 6 \sum_{y \geq 2} \sum_{\ell' \geq 0} S(y, \ell') a(x - y - 1, \ell - \ell') \right\}.$$

Indeed, if one of the first three cases holds, then

$$\mathbb{P}_p(CA(R^+) \geq x | W \subset A) \leq a(x - 2, \ell - 1),$$

by Lemma 22 (applied with $x = 2$) and the Crossing Algorithm. [Recall that $a(x, \ell)$ is monotone in both $x$ and $\ell$.] On the other hand, consider the fourth case, and recall that the event "there is a $y$-saver of an $L$-blocker at point $(1, j)$" means that the six sets $M_1(1)$, $M_1(2)$, $M_2(1)$, $M_1(j)$, $M_1(j - 1)$ and $M_2(j)$ are empty, and that at least one of them has a $y$-saver. Thus, our $y$-saver may lie either in $R[1, y]$ or in $R[2, y + 1]$, and so

$$\mathbb{P}_p(CA(R^+) \geq x | W \subset A, \exists y\text{-saver at point } (1, j)) \leq a(x - y - 1, \ell - \ell'),$$

where $\ell' = |W \cap R[1, y + 1]|$. The probability such a $y$-saver exists is at most $6S(y, \ell')$, so

$$\mathbb{P}_p(CA(R^+) \geq x | W \subset A) \leq 6 \sum_{y \geq 2} \sum_{\ell' \geq 0} S(y, \ell') a(x - y - 1, \ell - \ell'),$$

as required. Finally, in the fifth case $\mathbb{P}_p(CA(R^+) \geq x) = 0$.

Next suppose that $W \cap R[1, 2] = \varnothing$. Then $\mathbb{P}_p(CA(R^+) \geq x)$ is bounded above by $a(x - 1, \ell)$ in the first case, and by $a(x - 2, \ell)$ in the second and third cases, by Lemma 22. Moreover, the probability that the first case occurs is $2u - u^2$ and the probability of the second case is $u(1 - u)^2$.

Now recall that $1 - u = (1 - p)^m \geq e^{-2pm} \geq e^{-2B}$, since $m \leq B/p$. Thus, the probability that the third case occurs is at most

$$(1 - (1 - u)^3)^{(k-2)/2} \leq (1 - e^{-6B})^{(k-2)/2} \leq \delta^3,$$

since $k \geq 6e^{6B} \log(1/\delta) + 2$. The fourth and fifth cases are as before.  □

Next we use Lemma 26 and the induction hypothesis to bound the sum in the claim. Indeed, if $y \geq 2$ and $\ell' \geq 0$, then

$$S(y, \ell') a(x - y - 1, \ell - \ell')$$
$$\leq 4m(y+1)^3 (4kyp)^{\lceil (y+1)/2 \rceil - \ell'} (\beta(u) + \delta)^{x - y - 1 - 2(\ell - \ell')}.$$

Recall that $y \leq s \leq T/p$, and so $8kyp \leq \delta^4$ if $T$ is chosen to be sufficiently small compared with $k$ and $\delta$. Note also that $mk^2p^2 \leq Bk^2p \to 0$ as $p \to 0$. Thus,

$$\sum_{y \geq 2} \sum_{\ell' \geq 0} S(y, \ell') a(x - y - 1, \ell - \ell')$$
$$\leq 2 \sum_{y \geq 2} 4m(y+1)^3 (4kyp)^{\lceil (y+1)/2 \rceil} (\beta(u) + \delta)^{x - y - 1 - 2\ell}$$
$$\leq 2^{12} m(12kp)^2 (\beta(u) + \delta)^{x - 2\ell - 4} \leq \delta^3 (\beta(u) + \delta)^{x - 2\ell - 2}.$$

Finally, recalling that

$$\beta(u)^2 = (2u - u^2)\beta(u) + u(1 - u)^2,$$

and using the claim, the bounds above and the induction hypothesis, we get

$$\mathbb{P}_p(CA(R^+) \geq x | W \subset A)$$
$$\leq (\beta(u) + \delta)^{x - 2\ell - 2} ((2u - u^2)(\beta(u) + \delta) + u(1 - u)^2 + 7\delta^3)$$
$$\leq (\beta(u) + \delta)^{x - 2\ell - 2} (\beta(u)^2 + (2u - u^2)\delta + \delta^2)$$
$$\leq (\beta(u) + \delta)^{x - 2\ell},$$

since $\delta < \delta_0(B) \leq 1/7$ and $2u - u^2 < \beta(u)$.  □

Lemma 21 now follows easily from Lemmas 24–26.

Proof of Lemma 21. The lower bound follows easily by Lemma 7, applied with $d = 2$ and $\ell = 1$, since if $R$ is semi-crossed (in the sense of Section 3), then it is crossed (in the sense of this section). Note that $\dim(R) = (s, m)$, so $v(1) = m$ and $u(1) = u$, as required. We shall therefore concentrate on the upper bound.

Let $B > 0$ and $\delta_0(B) > \delta > 0$ be sufficiently small. We may assume that $\beta(u) + \delta < 1$, since otherwise the result is trivial. Let $k \geq 6e^{6B} \log(1/\delta) + 2$ and let $T$ be chosen appropriately so that Lemmas 25 and 26 hold.

Let $p > 0$ be sufficiently small, and let $R \subset C(n, k)$ be a rectangle as described, with $\dim(R) = (s, m)$, where $m \leq B/p$ and $s \leq T/p$. Let $2\ell \leq s$, let $W \subset R$ with $|W| \leq \ell$, and let $A \in \text{Bin}(R, p)$. We are required to show that

$$\mathbb{P}(R \text{ is left-to-right crossed by } A | W \subset A) \leq m^3 (\beta(u) + \delta)^{s - 4\ell}.$$

Indeed, suppose that $R$ is crossed by $A$. Then, by Lemma 24, there exists an $x \in [s]$ such that the events,

(a) $R$ is unblocked up to the point $x$, and
(b) there is a last chance at distance $y = s - x$,

occur disjointly. Therefore, by Reimer's Theorem,

$$\mathbb{P}(R \text{ is left-to-right crossed by } A | W \subset A) \le \sum_{x=0}^{s} a(x, \ell) b(s - x, \ell).$$

Letting $y = s - x$, and applying Lemmas 25 and 26, we get

$$a(x, \ell) b(y, \ell) \le (\beta(u) + \delta)^{x - 2\ell} \min\{2my(6kyp)^{y/2 - \ell}, 1\} \le m^2(\beta(u) + \delta)^{s - 4\ell},$$

since we may choose $T$ small enough that $6kyp \le 6kT \le \delta^2$. Indeed, either $y < 2\ell$, in which case $x - 2\ell > s - 4\ell$, or $y \ge 2\ell$, in which case

$$2my(6kyp)^{y/2 - \ell} \le m^2 \delta^{y - 2\ell} \le m^2(\beta(u) + \delta)^{y - 2\ell}.$$

Hence,

$$\mathbb{P}(R \text{ is left-to-right crossed by } A | W \subset A) \le m^3(\beta(u) + \delta)^{s - 4\ell},$$

as required.  $\square$

4.3. *Proof of Theorem 13.*  In this section we complete the proof of Theorem 13, using Lemmas 20 and 21 and the method of Holroyd [19]. We begin by bounding the probability that a rectangle grows sideways by $T/p$.

Let $R \subset R'$ be rectangles in $C(B/p, k)$, and recall from Section 4.1 the definition of $D(R, R')$, and that

$$P_p(R, R') = \mathbb{P}(D(R, R') | A \in \text{Bin}(R', p)).$$

Let $R_1, \dots, R_8$ be as in Figure 2. Moreover, let $R_{\text{top}} = R_1 \cup R_2 \cup R_3$, $R_{\text{right}} = R_3 \cup R_4 \cup R_5$, $R_{\text{bottom}} = R_5 \cup R_6 \cup R_7$ and $R_{\text{left}} = R_1 \cup R_7 \cup R_8$. We have

$$P_p(R, R') = \mathbb{P}_p(R' \text{ is internally spanned by } A | R \text{ is internally filled})$$

$$\le \mathbb{P}_p(H^{\rightarrow}(R_{\text{right}}) \cup H^{\leftarrow}(R_{\text{left}}) \cup H^{\uparrow}(R_{\text{top}}) \cup H^{\downarrow}(R_{\text{bottom}})).$$

Let $g(z) = g_2(z)$, the function defined in the Introduction. We shall deduce the following lemma from Lemma 21.

Lemma 28.  *Let $B > 0$, and let $\delta_0(B) > \delta > 0$, $Z_0 > Z > 0$ and $k \in \mathbb{N}$ satisfy $k \ge 12 e^{6B} \log(1/\delta) + 2$ and $(6\delta)^2 < Z$. Then there exists a constant $T = T(Z, k, \delta) > 0$ such that the following holds.*

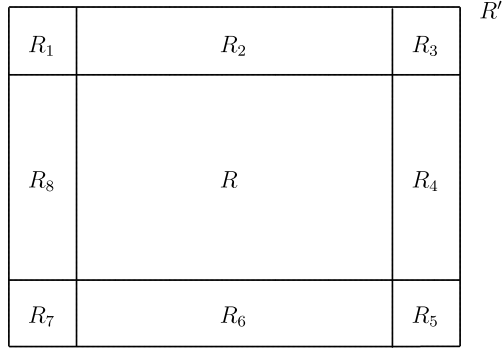

FIG. 2.    *The rectangles* $R \subset R'$.

Let $p > 0$ be sufficiently small, and let $R \subset R' \subset C(B/p, k)$ be rectangles, with $\dim(R) = (m, n)$ and $\dim(R') = (m+s, n+t)$, where $Z/p \leq m, n \leq B/p$ and $s, t \leq T/p$. Then

$$(6) \quad \begin{aligned} P_p(R, R') &\leq m^7 n^7 (\beta(u(n+t)) + \delta^2)^{(1-\delta)s} (\beta(u(m+s)) + \delta^2)^{(1-\delta)t} \\ &\leq (B/p)^{14} \exp(-(1-2\delta)(g(qn)s + g(qm)t)), \end{aligned}$$

where $u(x) = 1 - (1-p)^x = 1 - e^{-qx}$.

PROOF.    Let $C = R_1 \cup R_3 \cup R_5 \cup R_7$ denote the corner areas of $R' \setminus R$, and let $W = A \cap C$. Let $\ell = |W|$, and note that $|C| = st \leq (T/p)^2$. The idea is that, since $T$ may be chosen small compared with $Z$, it is likely that $\ell$ will be small compared with $s$ and $t$, and so the events $H^{\rightarrow}(R_{\text{right}})$, $H^{\leftarrow}(R_{\text{left}})$, $H^{\uparrow}(R_{\text{top}})$ and $H^{\downarrow}(R_{\text{bottom}})$ are "almost independent."

To be precise, let us apply Lemma 21 to (appropriate rotations of) the rectangles $R_{\text{right}}$, $R_{\text{left}}$, $R_{\text{top}}$ and $R_{\text{bottom}}$, conditional on the event that $|A \cap C| = \ell$. Let $u_1 = u(n+t)$ and $u_2 = u(m+s)$, and assume that $\delta$ is sufficiently small, so, in particular, $\beta(u_i) + \delta^2 < 1$ for $i = 1, 2$. Let $T = T(k, \delta^2)$ be chosen small enough so that Lemma 21 holds, and so that $(m+s)^6(n+t)^6 \leq 2m^6 n^6$. Then, by Lemma 21, applied to $B$, $\delta^2$ and $k$,

$$\begin{aligned} \mathbb{P}_p(D(R, R') || W| = \ell) &\leq h(R_{\text{right}}, \ell) h(R_{\text{left}}, \ell) h(R_{\text{top}}, \ell) h(R_{\text{bottom}}, \ell) \\ &\leq 2m^6 n^6 (\beta(u_1) + \delta^2)^{s-8\ell} (\beta(u_2) + \delta^2)^{t-8\ell}. \end{aligned}$$

We split into two cases: $8\ell < \delta \min\{s, t\}$ and $8\ell \geq \delta \min\{s, t\}$. In the first case we have

$$(7) \quad \begin{aligned} &\mathbb{P}_p(D(R, R') \text{ and } 8|W| < \delta \min\{s, t\}) \\ &= \sum_{8\ell < \delta \min\{s,t\}} \mathbb{P}_p(D(R, R') || W| \leq \ell) \end{aligned}$$

$$\leq \delta(s+t)m^6 n^6 (\beta(u_1) + \delta^2)^{s-8\ell}(\beta(u_2) + \delta^2)^{t-8\ell}$$

$$\leq \frac{m^7 n^7}{2}(\beta(u(n+t)) + \delta^2)^{(1-\delta)s}(\beta(u(m+s)) + \delta^2)^{(1-\delta)t}.$$

In the second case, first note that $|W| \sim \text{Bin}(st, p)$, and so, if $8\ell \geq \delta \min\{s, t\}$, then

$$\mathbb{P}_p(|W| = \ell) \leq \binom{st}{\ell}p^\ell \leq \left(\frac{3pst}{\ell}\right)^\ell \leq \left(\frac{24pst}{\delta \min\{s, t\}}\right)^\ell \leq \left(\frac{50T}{\delta}\right)^\ell.$$

Thus,

$$\begin{aligned}
&\mathbb{P}_p(D(R, R') \text{ and } 8|W| \geq \delta \min\{s, t\}) \\
&\qquad = \sum_{8\ell \geq \delta \min\{s,t\}} \mathbb{P}_p(D(R, R')||W| = \ell)\mathbb{P}_p(|W| = \ell) \\
&\qquad \leq \sum_{8\ell \geq \delta \min\{s,t\}} 2m^6 n^6 (\beta(u_1) + \delta^2)^{s-8\ell}(\beta(u_2) + \delta^2)^{t-8\ell}\left(\frac{50T}{\delta}\right)^\ell \\
&\qquad \leq 4m^6 n^6 (\beta(u(n+t)) + \delta^2)^s (\beta(u(m+s)) + \delta^2)^t,
\end{aligned}$$

(8)

since $T$ may be chosen so that $100T \leq \delta^{33}$. The first inequality follows from (7) and (8).

To obtain the second inequality, recall that $m, n \leq B/p$, and that

$$e^{-g(qx)} = \beta(1 - e^{-qx}) = \beta(u(x)).$$

Now, recall that $m, n \geq Z/p$, so $u(m+s), u(n+t) \geq Z/2$, that $\beta(u)$ is increasing in $u$, and that $\beta(u) \geq \sqrt{u}/2$ for small $u$. Thus, $\beta(u_i) \geq \sqrt{Z}/3$, and so, since $(6\delta)^2 < Z$,

$$\beta(u_i) + \delta^2 \leq \beta(u_i)^{1-\delta}$$

for $i = 1, 2$. [Note that $\frac{d}{dx}(c + x^2 - c^{1-x}) = 2x + c^{1-x}\log c < 0$ if $0 < 2x < c < 1/e$.] Therefore,

$$(\beta(u_1) + \delta^2)^{(1-\delta)s} \leq \beta(u_1)^{(1-\delta)^2 s} = e^{-(1-\delta)^2 sg(q(n+t))} \leq \exp(-(1-2\delta)sg(qn))$$

since $g(x)$ is continuous, and $qt \leq 2T$ may be made arbitrarily small compared with $qn \geq Z$ and $\delta$. A similar inequality holds for $\beta(u_2)$ and $t$, and so the result follows. $\square$

We now rewrite the right-hand side of (6) in a more useful form. Define

$$W_g(\mathbf{a}, \mathbf{b}) = \inf_{\gamma: \mathbf{a} \to \mathbf{b}} \int_\gamma (g(y)\, dx + g(x)\, dy),$$

where the infimum is taken over all piecewise linear, increasing paths from $\mathbf{a}$ to $\mathbf{b}$ in $\mathbb{R}^2$ (see Section 6 of [19]). Moreover, for any two rectangles $R \subset R'$, let

$$U(R, R') = W_g(q \dim(R), q \dim(R')).$$

The following easy observation holds not only for $g$, but for any decreasing function.

OBSERVATION 29 (Proposition 12 of [19]).

$$W_g(\mathbf{a}, \mathbf{b}) \geq (b_1 - a_1)g(a_2) + (b_2 - a_2)g(a_1).$$

The following corollary of Lemma 28 is now immediate.

COROLLARY 30. *Under the conditions of Lemma* 28,

$$P_p(R, R') \leq (B/p)^{14} \exp\left(-\frac{(1 - 2\delta)U(R, R')}{q}\right).$$

Next we bound the probability that a seed is internally spanned.

LEMMA 31. *Let* $\alpha > 0$, $Z > 0$ *and* $k \in \mathbb{N}$, *with* $2kZ \leq e^{-4\alpha}$. *Let* $n \in \mathbb{N}$ *and* $p > 0$, *let* $R \subset C(n, k)$ *be a rectangle with* $\mathrm{short}(R) \leq Z/p$, *and let* $A \in \mathrm{Bin}(R, p)$. *Then*

$$\mathbb{P}(R \in \langle A \rangle) \leq e^{-\alpha \phi(R)}.$$

PROOF. Suppose $\dim(R) = (u, v)$, with $u \leq v$. Note that if $R \in \langle A \rangle$, then $R$ has no double gap. Thus, by Lemma 27,

$$\mathbb{P}(R \in \langle A \rangle) \leq (2kup)^{v/2} \leq (2kZ)^{\phi(R)/4} \leq e^{-\alpha \phi(R)},$$

as required. □

Finally, in order to deduce Theorem 13 from Corollary 30 and Lemmas 19, 20 and 31, we shall need some way to relate the quantities $\sum_{\vec{\Gamma}(u) = \{v\}} U(R_v, R_u)$ and $\sum_{\text{seeds } u} \phi(R_u)$. The following lemma, due to Holroyd [19], does this for us.

LEMMA 32 (Lemma 37 of [19]). *Let* $n, k \in \mathbb{N}$ *and* $T, Z, p > 0$. *For any hierarchy* $\mathcal{H}$ *of a rectangle* $R \subset C(n, k)$ *which is good for* $(T, Z, p)$, *there exists a rectangle* $S = S(\mathcal{H}) \subset R$, *with*

$$\phi(S) \leq \sum_{\text{seeds } u} \phi(R_u)$$

*such that*

$$\sum_{\vec{\Gamma}(u)=\{v\}} U(R_v, R_u) \geq U(S, R) - (2qg(Z))|\{u \in \mathcal{H} : |\vec{\Gamma}(u)| \geq 2\}|.$$

PROOF. This was proved in [19] only for $[n]^2$, but the proof for $C(n, k)$ is exactly the same.  □

Finally, we need the following simple modification of a lemma from [19].

LEMMA 33 (Proposition 14 of [19]). *If $a_1 + a_2 \leq A$ and $\mathbf{b} = (B, b_2)$, and $a_2 \leq b_2$, then*

$$W(\mathbf{a}, \mathbf{b}) \geq 2\int_A^B g(z)\, dz - Bg(B).$$

We remark that $Bg(B) \to 0$ as $B \to \infty$, since $g$ is integrable on $(0, \infty)$, or by Proposition 3.

We are finally ready to prove Theorem 13.

PROOF OF THEOREM 13. Let $\varepsilon > 0$, and let $B = B(\varepsilon)$, $\delta$, $k$, $\alpha$, $Z$ and $T = T(Z, k, \delta)$ be positive constants, chosen so that Lemmas 21, 28 and 31 all hold. Thus, $B$, $k$ and $\alpha$ are sufficiently large, and $\delta$, $Z$ and $T$ are sufficiently small. In particular, let $\alpha = 2B$, $k \geq 10e^{6B}\log(1/\delta)$, $6\delta^2 \leq Z$ and $kZ \leq e^{-4\alpha}$. It is easy to see that we can satisfy these inequalities simultaneously, and that we have

$$T \ll \delta, Z \ll 1 \ll B \ll k.$$

Finally, we let $p \to 0$, so $p \ll T$.

Let $R \subset C(B/p, k)$ with $\text{long}(R) = B/p$, and let $A \in \text{Bin}(R, p)$. By Corollary 30 and Lemmas 19, 20 and 31, we obtain

$$\mathbb{P}(R \in \langle A \rangle)$$

$$\leq \sum_{\mathcal{H} \in \mathcal{H}(R,T,Z,p)} \left( \prod_{\vec{\Gamma}(u)=\{v\}} P_p(R_v, R_u) \right) \prod_{\text{seeds } u} P_p(R_u)$$

$$\leq Mp^{-M}(B/p)^{14M}$$

$$\times \exp\left( -\sum_{\vec{\Gamma}(u)=\{v\}} \frac{(1-2\delta)U(R_v, R_u)}{q} - \alpha \sum_{\text{seeds } u} \phi(R_u) \right).$$

Now, applying Lemma 32, this becomes

$$\mathbb{P}(R \in \langle A \rangle) \leq M'p^{-M'}\exp\left( -\frac{(1-2\delta)U(S, R)}{q} - \alpha \sum_{\text{seeds } u} \phi(R_u) \right)$$

for some constant $M'$, since $g(Z)|V(G_{\mathcal{H}})|$ is bounded above by a constant depending only on $B$, $Z$ and $T$.

We split into two cases, depending on whether $\sum_{\text{seeds } u} \phi(R_u) \geq \frac{1}{Bp}$ or not, and apply Lemma 33. In the former case, we get

$$\mathbb{P}(R \in \langle A \rangle) \leq M'p^{-M'} \exp\left(-\frac{(1-2\delta)U(S,R)}{q} - \alpha \sum_{\text{seeds } u} \phi(R_u)\right)$$

$$\leq M'p^{-M'} \exp\left(-\frac{\alpha}{Bp}\right) \leq \exp\left(-\frac{2\lambda(3,3)}{p}\right)$$

since $\alpha = 2B$, and $p$ is sufficiently small. In the latter case note that $\phi(S) \leq \frac{1}{Bp} \leq \frac{2}{Bq}$, and that $xg(x) \to 0$ as $x \to \infty$. Thus, by Lemma 33,

$$\mathbb{P}(R \in \langle A \rangle) \leq M'p^{-M'} \exp\left(-\frac{(1-2\delta)U(S,R)}{q} - \alpha \sum_{\text{seeds } u} \phi(R_u)\right)$$

$$\leq M'p^{-M'} \exp\left(-\frac{2(1-2\delta)}{q}\left(\int_{2/B}^{B} g(z)\,dz - Bg(B)\right)\right)$$

$$\leq \exp\left(-\frac{2\lambda(3,3)-\varepsilon}{p}\right)$$

if $B$ is sufficiently large and $\delta$ and $p$ are sufficiently small, as required. $\quad\square$

4.4. *Proofs of Corollary 14 and Theorem 2.* We complete Section 4 by making the easy final steps necessary to deduce Corollary 14 and Theorem 2.

PROOF OF COROLLARY 14. Let $\varepsilon > 0$, and let $B = B(\varepsilon)$, $k = k(B,\varepsilon)$ be chosen according to Theorem 13. Let $n = n(B,k,\varepsilon)$ be sufficiently large, let $p = \frac{\lambda(3,3)-\varepsilon}{\log n}$, and let $A \in \text{Bin}(C(n,k),p)$. We are required to show that

$$\mathbb{P}(\text{long}(R) \geq B\log n \text{ for some } R \in \langle A \rangle) \leq n^{-\varepsilon}.$$

Indeed, suppose $\text{long}(R) \geq B\log n$ for some $R \in \langle A \rangle$. By Lemma 16, there exists an internally spanned rectangle $R' \subset R$ with $(B/2)\log n \leq \text{long}(R') \leq B\log n$. Then, by Theorem 13,

$$\mathbb{P}(R' \in \langle A \cap R' \rangle) \leq \exp\left(-\frac{2\lambda(3,3)-\varepsilon}{p}\right)$$

$$= \exp\left(-\left(\frac{2\lambda(3,3)-\varepsilon}{\lambda(3,3)-\varepsilon}\right)\log n\right) \leq n^{-(2+2\varepsilon)},$$

if $B$ is sufficiently large, since $\lambda(3,3) < 1/2$.

There are at most $(B\log n)^2 n^2 \leq n^{2+\varepsilon}$ potential such rectangles $R'$. So, writing $X(B)$ for the number of internally spanned rectangles $R' \subset C(n,k)$ with $(B/2)\log n \leq \text{long}(R') \leq B\log n$, we get

$$\mathbb{P}(\text{long}(R) \geq B\log n \text{ for some } R \in \langle A \rangle) \leq \mathbb{E}(X(B)) \leq n^{-\varepsilon},$$

as required. $\square$

The proof of Theorem 2 is now immediate.

PROOF OF THEOREM 2. Let $\varepsilon > 0$, let $k \geq K(\varepsilon) \in \mathbb{N}$ be chosen according to Corollary 14, and let $\delta > 0$ satisfy $\delta \leq (1-\delta)^k$. Finally, let $n \in \mathbb{N}$ be sufficiently large. We are required to prove that

$$\frac{\lambda(3,3) - \varepsilon}{\log n} \leq p_\delta(C(n,k)) \leq p_{1-\delta}^{(s)}(C^*(n,2)) \leq \frac{\lambda(3,3) + \varepsilon}{\log n}.$$

The lower bound is immediate from Corollary 14, which says that, moreover, if $p = \frac{\lambda(3,3)-\varepsilon}{\log n}$, then with high probability (as $n \to \infty$) $[A]$ has diameter only $O(\log n)$.

The middle bound follows from the condition $\delta \leq (1-\delta)^k$, and the fact that there exist two copies ($C_1$ and $C_2$, say) of $C^*(n,2)$ in $C(n,k)$. Indeed, if $p = p_{1-\delta}^{(s)}(C^*(n,2))$, then $\mathbb{P}_p(A$ percolates in $C(n,k)) \geq (1-\delta)^k$ for any $2 \leq k \in \mathbb{N}$.

To spell it out, we use induction on $k$. The result holds for the base case, $k = 2$, because all sites in $C(n,2)$ have threshold 2, so we may couple $C(n,2)$ with two overlapping copies of $C^*(n,2)$, and use the FKG inequality.

For the induction step, let $C_1$ be the copy of $C^*(n,2)$ in $C(n,k)$ with vertex set $\{(x,y,z) : z \in \{1,2\}\}$, and let $D$ be a copy of $C(n,k-1)$ on vertex set $\{(x,y,z) : 2 \leq z \leq k\}$. Observe that if $A \cap C_1$ semi-percolates in $C_1$, and $A \cap D$ percolates in $D$, then $A$ percolates in $C(n,k)$. Moreover, these events are increasing in $A$, and have probability at least $1-\delta$ and $(1-\delta)^{k-1}$ respectively, by the induction hypothesis. Thus, the result follows by the FKG inequality.

Finally, the upper bound was proved in Section 3. Indeed, Theorem 5 (applied in the case $d = r = 2$, $\ell = 1$) says exactly that, if $p = \frac{\lambda(3,3)+\varepsilon}{\log n}$, then

$$\mathbb{P}_p(A \subset C^*(n,2) \text{ semi-percolates}) \to 1$$

as $n \to \infty$. $\square$

## 5. Proof of Theorem 1.

In this section we shall use Corollary 14 to prove Theorem 1. To do so, we will borrow the ideas of Cerf, Cirillo and Manzo [12, 13], and also of Holroyd [20], who corrects a small error from [12, 13].

In order to state the main lemma of this section, we need a little notation. We will be interested in two-colored graphs, that is, simple graphs with two types of edges, which we shall label "good" and "bad." We call such a two-colored graph "admissible" if it either contains at least one bad edge, or if every component is a clique. For any set $S$, let

$$\Lambda(S) := \{\text{admissible two-colored graphs with vertex set } S \times [2]\}.$$

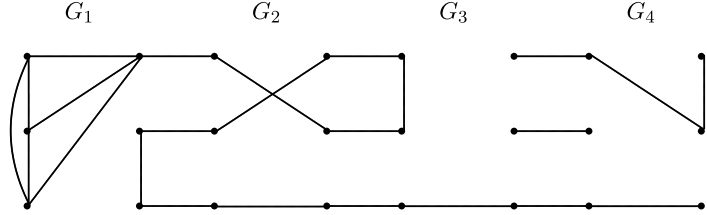

$G_1$          $G_2$          $G_3$          $G_4$

Fig. 3.  *A graph $G_{\mathcal{P}}$, with $S = [3]$ and $m = 4$.*

Now, given $m \in \mathbb{N}$, let

$$\Omega(S, m) := \{\mathcal{P} = (G_1, \ldots, G_m) : G_t \in \Lambda(S) \text{ for each } t \in [m]\},$$

the set of sequences of two-colored admissible graphs on $S \times [2]$ of length $m$. We shall sometimes think of $G_t$ as a colored graph on $S \times [2t-1, 2t]$, and trust that this will cause no confusion. We shall be interested in probability distributions on $\Omega(S, m)$ in which, with high probability, there are bad edges in only very few of the graphs $G_t$.

Now, for each $\mathcal{P} \in \Omega(S, m)$, let $G_{\mathcal{P}}$ denote the graph with vertex set $S \times [2m]$, and the following edge set $E(G_{\mathcal{P}})$ (see, e.g., Figure 3):

(a)  $G_{\mathcal{P}}[S \times \{2y-1, 2y\}] = G_y$,
(b)  $\{(x, 2y), (x', 2y+1)\} \in E(G_{\mathcal{P}}) \Leftrightarrow x = x'$,
(c)  $\{(x, y), (x', y')\} \notin E(G_{\mathcal{P}})$ if $|y - y'| \geq 2$.

Edges in $G_{\mathcal{P}}$ of types (a) and (b) are labeled good and bad in the obvious way, to match the label of the corresponding edge in $G_y$. Thus, $G_{\mathcal{P}}$ has three types of edges: good, bad and unlabeled.

Such a graph $G_{\mathcal{P}}$, with $S = [3]$ and $m = 4$, is pictured below. Note that, for example, $G_3$ has two edges: $\{(1,1), (2,1)\}$ and $\{(3,1), (3,2)\}$, and that $G_4$ must contain a bad edge.

Given $G \in \Lambda(S)$, let $E^g(G)$ denote the set of good edges, and $E^b(G)$ denote the bad edges, so that $E(G) = E^g(G) \cup E^b(G)$. If $uv$ is a good edge in $G$, then we shall write $u \sim v$.

For each $1 \leq i \leq j \leq 2m$, we shall write $G_{\mathcal{P}}[i, j]$ for the subgraph of $G_{\mathcal{P}}$ induced by the set $S \times [i, j]$, and $V(G_t)$ for the vertex set of $G_{\mathcal{P}}[2t-1, 2t]$. For each vertex $v = (x, y) \in V(G_{\mathcal{P}})$, let

$$\Gamma_{\mathcal{P}}(v) := \{u \in V(G_{\lceil y/2 \rceil}) : u \sim v \text{ and } u \neq v\},$$

and let $d_{\mathcal{P}}(v) = |\Gamma_{\mathcal{P}}(v)|$. We emphasize that $d_{\mathcal{P}}(v)$ is the number of *good* edges incident with $v$.

We shall use the following simple calculation in the proof below.

LEMMA 34. *Let $m, r \in \mathbb{N}$ and, for each $t \in [r]$, let $i_t, j_t \in \{1, 2\}$ and $k(t) \in [m]$. Let $S$ be any finite set, and $\mathcal{P} = (G_1, \ldots, G_m)$ be a random sequence of admissible two-colored graphs on $S \times [2]$, chosen according to some (arbitrary) probability distribution $f_\Omega$ on $\Omega(S, m)$. Then*

$$\sum_{x_1, \ldots, x_{r+1} \in S} \prod_{t=1}^{r} \mathbb{P}((x_t, i_t) \sim (x_{t+1}, j_t) \ in \ G_{k(t)}) \leq |S| \Big( \max_{v \in G_\mathcal{P}} \mathbb{E}(d_\mathcal{P}(v)) \Big)^r.$$

PROOF. This follows easily from the fact that $\mathbb{E}(d_\mathcal{P}(v)) = \sum_u \mathbb{P}(u \sim v)$. Indeed, pulling constant factors through the summation signs, the left-hand side may be rewritten as

$$\sum_{x_1, x_2} \Big( \mathbb{P}((x_1, i_1) \sim (x_2, j_1)) \sum_{x_3} \Big( \mathbb{P}((x_2, i_2) \sim (x_3, j_2)) \cdots$$
$$\times \sum_{x_{r+1}} \mathbb{P}((x_r, i_r) \sim (x_{r+1}, j_r)) \Big) \ldots \Big),$$

where $|S|$ is the number of choices for $x_1$. The result now follows by using the inequalities

$$\sum_{x_{t+1}} \mathbb{P}((x_t, i_t) \sim (x_{t+1}, j_t) \ in \ G_{k(t)}) \leq \max_{v \in G_\mathcal{P}} \mathbb{E}(d_\mathcal{P}(v))$$

for each $t \in [r]$.  □

Finally, let $X(\mathcal{P})$ denote the event that there is a connected path across $G_\mathcal{P}$ (i.e., a path from the set $S \times \{1\}$ to the set $S \times \{2m\}$). Observe that the event $X(\mathcal{P})$ holds for the graph $G_\mathcal{P}$ depicted in Figure 3.

The following lemma is proved, but not stated, by Cerf and Cirillo [12], and by Cerf and Manzo [13]. Since it is not immediately obvious how to read the result out of their papers, we give a sketch of the proof.

LEMMA 35 (Cerf and Cirillo [12]). *For each $0 < \alpha < 1/2$ and $\varepsilon > 0$, there exists $\delta > 0$ such that the following holds for all $m \in \mathbb{N}$ and all finite sets $S$ with $\alpha^4 |S|^\varepsilon \geq 1$.*

*Let $\mathcal{P} = (G_1, \ldots, G_m)$ be a random sequence of admissible two-colored graphs on $S \times [2]$, chosen according to some probability distribution $f_\Omega$ on $\Omega(S, m)$. Suppose $f_\Omega$ satisfies the following conditions:*

(a) *independence: $G_i$ and $G_j$ are independent if $i \neq j$,*
(b) *BK condition: for each $t \in [m]$, $r \in \mathbb{N}$, and each $x_1, y_1, \ldots, x_r, y_r \in V(G_t)$,*

$$\mathbb{P}\Big( \bigwedge_{j=1}^{r} (x_j \sim y_j) \wedge \bigwedge_{j \neq j'} (x_j \not\sim x_{j'}) \wedge (E^b(G_t) = \varnothing) \Big)$$
$$\leq \prod_{j=1}^{r} \mathbb{P}(x_j \sim y_j),$$

*and for each $t \in [m]$ and $v \in V(G_\mathcal{P})$,*

(c) *bad edge condition:* $\mathbb{P}(E^b(G_t) \neq \varnothing) \leq |S|^{-\varepsilon}$,

(d) *good edge condition:* $\mathbb{E}(d_\mathcal{P}(v)) \leq \delta$.

*Then*

$$\mathbb{P}(X(\mathcal{P})) \leq \alpha^m |S|.$$

REMARK 7. In our application $S$ will be the set $[n]^2$, and will correspond to the top (or bottom) layer of a copy of $C(n,k)$. The pair $uv$ will be an edge of the graph $G_t$ if $u,v \in C(n,k)$ are in the same component of $[A]$, where $A \in \text{Bin}(C(n,k),p)$. Edges will be labeled "good" if both endpoints lie in some internally filled component of "small" diameter, that is, less than $B \log n$, where $B > 0$ is sufficiently large.

Condition (b) will be proved using the van den Berg–Kesten Lemma, condition (c) using Corollary 14, and condition (d) by Lemma 36, below.

PROOF OF LEMMA 35. Let $\{z_1, \ldots, z_t\} \subset [m]$ denote the indices for which $E^b(G_z) \neq \varnothing$, and note that this event has probability at most $n^{-\varepsilon t}$, where $n := |S|$. Thus, the probability that $t \geq T := 3 \log(1/\alpha)m/(\varepsilon \log n)$ is at most

$$2^m n^{-3 \log(1/\alpha)m/\log n} \leq \alpha^{m+1}.$$

So suppose $t \leq T$; for each pair $z_j$, $z_{j+1}$, we shall count "shortest" paths between the left- and right-hand sides of $G_\mathcal{P}[1,2s] \cong G_\mathcal{P}[2z_j + 1, 2(z_{j+1}-1)]$.

Indeed, let $\tilde{X}(2s)$ denote the event that there is a path across $G_\mathcal{P}[1,2s]$, *and* that $E^b(G_z) = \varnothing$ for each $z \in [s]$. We claim that if $\tilde{X}(2s)$ holds, then there is a sequence of (distinct) vertices $(x_1, i_1), (y_1, j_1), \ldots, (x_r, i_r), (y_r, j_r) \in S \times [1, 2s]$, with $r \geq s$, such that:

- $i_1 = 1$ and $j_r = 2s$,
- $x_{t+1} = y_t$ and

$$i_{t+1} - j_t = \begin{cases} 1, & \text{if } j_t \equiv 0 \pmod 2, \\ -1, & \text{if } j_t \equiv 1 \pmod 2, \end{cases} \qquad \text{for each } t \in [r-1],$$

- $(x_t, i_t) \sim (y_t, j_t)$ for each $t \in [r]$,
- $(x_t, i_t) \not\sim (x_{t'}, i_{t'})$ for each $t \neq t'$.

Indeed, to obtain such a path for which these events occur, simply choose a path with $r$ minimal. (Note that we use here that each graph $G_t$ is admissible.) Let $\mathcal{J}$ denote the collection of such sequences, that is, the collection of minimal paths across $G_\mathcal{P}$. Summing over all sequences in $\mathcal{J}$, we get

$$\mathbb{P}(\tilde{X}(2s)) \leq \sum_\mathcal{J} \mathbb{P}\left( \bigwedge_{t=1}^r ((x_t, i_t) \sim (y_t, j_t)) \right)$$

$$\wedge \bigwedge_{t \neq t'} ((x_t, i_t) \not\sim (x_{t'}, i_{t'})) \wedge \bigwedge_{t=1}^{s} (E^b(G_t) = \varnothing) \Bigg).$$

Now, using conditions (a), (b) and (d), and Lemma 34, it follows that

$$\mathbb{P}(\tilde{X}(2s)) \leq \sum_{\mathcal{J}} \mathbb{P}\Bigg( \bigwedge_{t=1}^{r} ((x_t, i_t) \sim (y_t, j_t))$$

$$\wedge \bigwedge_{t \neq t'} ((x_t, i_t) \not\sim (x_{t'}, i_{t'})) \wedge \bigwedge_{t=1}^{s} (E^b(G_t) = \varnothing) \Bigg)$$

$$\leq \sum_{\mathcal{J}} \prod_{t=1}^{r} \mathbb{P}((x_t, i_t) \sim (y_t, j_t))$$

$$\leq \sum_{r \geq s} 4^r |S| \Big( \max_{v \in G_{\mathcal{P}}[1, 2s]} \mathbb{E}(d_{\mathcal{P}}(v)) \Big)^r \leq 2n(4\delta)^s,$$

assuming $\delta$ is sufficiently small. The term $4^r$ comes from summing over all choices of $i_1, j_1, \ldots, i_r, j_r$.

Now, we simply sum over all choices of the set $\{z_1, \ldots, z_t\}$. Recalling that $t \leq T = 3 \log(1/\alpha) m / (\varepsilon \log n) \leq 3m/4$, and writing $s(j) = z_{j+1} - z_j - 1$ for each $j \in [0, t]$ (let $z_0 = 0$ and $z_{t+1} = m + 1$), this gives

$$\mathbb{P}(X(\mathcal{P})) \leq 2^m \prod_{j=0}^{t} 2n(4\delta)^{s(j)} \leq 2^{4m+1} n^{1+3 \log(1/\alpha) m / (\varepsilon \log n)} \delta^{m-t}$$

$$\leq n \alpha^{-3m/\varepsilon} \delta^{m/4} \leq \alpha^{m+1} n,$$

as required, since we may choose $\delta = \delta(\alpha, \varepsilon)$ as small as we like.  $\square$

In order to apply Lemma 35, we need to give an upper bound on the expected number of good edges incident to any given vertex. The next lemma does this. Given a bootstrap structure $G$ on $[n]^d \times [k]^\ell$, a set $A \subset V(G)$, a vertex $x \in V(G)$ and a number $R > 0$, define

$$\Gamma_G(A, R, x) := \{y \in V(G) \colon \text{there exists an internally filled}$$

$$\text{connected component } X \subset V(G) \text{ such that}$$

$$x, y \in X \text{ and } \operatorname{diam}(X) \leq R\}.$$

(This definition is important, and is due to Holroyd [20].) The following lemma, together with Corollary 14, allows us to apply Lemma 35. Since the proof is the same, and we shall need the result in [7], we prove it in $C([n]^d \times [k]^\ell, 2)$.

LEMMA 36. *Let $2 \le n, d \in \mathbb{N}$ and $\ell \in \mathbb{N}_0$. There exists a function $f(B, k) = f_{d,\ell}(B, k)$ such that, for any $B > 0$, any $k \in \mathbb{N}_0$ and any sufficiently small $p > 0$, the following holds. Let $G = C([n]^d \times [k]^\ell, 2)$, $A \in \mathrm{Bin}(V(G), p)$, $R = B/p^{1/(d-1)}$, and $x \in V(G)$. Then*

$$\mathbb{E}(|\Gamma_G(A, R, x)|) \le f(B, k)(\log(1/p))^{3d+\ell+1} p.$$

PROOF. Let $A \in \mathrm{Bin}([n]^d \times [k]^\ell, p)$, let $x, y \in [n]^d \times [k]^\ell$, and suppose $\|x - y\|_\infty = m$. For each $t \in \mathbb{N}$, let $a(m, t)$ denote the maximal probability (over all such choices of $x$ and $y$) that there exists an internally filled connected component $X \subset [n]^d \times [k]^\ell$ such that $x, y \in X$ and $m + 1 \le t = \mathrm{diam}(X) \le B/p^{1/(d-1)}$. We shall bound $a(m, t)$ from above.

Indeed, suppose such a component $X$ exists, and let $t := \mathrm{diam}(X)$, so $m + 1 \le t \le B/p^{1/(d-1)}$. We claim that

$$a(m, t) \le t^{2d} k^{2\ell} (2t^d k^\ell p)^{\lceil (t+1)/2 \rceil}.$$

This bound follows by considering the smallest cuboid containing $X$. It has diameter $t$, it contains $x$, and, since it is the smallest cuboid containing an internally filled component, it has no double gaps. There are at most $t^{2d} k^{2\ell}$ cuboids of diameter $t$ containing $x$, and the probability each has no double gap is at most $(2t^d k^\ell p)^{\lceil (t+1)/2 \rceil}$. (This follows exactly as in the proof of Lemma 27.)

The above bound works well for small $t$; for larger $t$ we use the bound

$$a(m, t) \le dt \exp(-\delta t)$$

for some $\delta = \delta(B, k) > 0$, which follows because $t \le B/p^{1/(d-1)}$. Indeed, let $t$ be as above and choose vertices $u, v \in X$ with $\|u - v\|_\infty + 1 = t$. Assume $t \ge k$, so that, without loss of generality, $u$ and $v$ differ by $t - 1$ in direction 1. Then the cuboid with dimensions $[t] \times [2B/p^{1/(d-1)}]^{d-1} \times [k]^\ell$, centered on $x$, with $u$ contained in one face and $v$ in the opposite face, has no double gap in direction 1. There are $dt$ such cuboids, and so the probability that there exists such a cuboid with no double gap is at most

$$dt(1 - (1 - p)^{2^d B^{d-1} k^\ell / p})^{t/2} \le dt(1 - \delta)^t \le dt \exp(-\delta t),$$

if $\delta = \delta(B, k) > 0$ is chosen to be sufficiently small, as claimed.

Now, there are at most $(2m + 1)^{d+\ell}$ vertices at (infinity norm) distance exactly $m$ from $x$, and, hence,

$$\mathbb{E}(|\Gamma_G(A, R, x)|) \le \sum_{m=0}^{R} (2m + 1)^{d+\ell} \sum_{t=m+1}^{R} a(m, t),$$

where $R = B/p^{1/(d-1)}$. Let $M = \frac{5(d+\ell)\log(1/p)}{\delta}$. The first bound on $a(m,t)$ gives

$$\sum_{t=m+1}^{M} a(m,t) \leq \sum_{t=m+1}^{M} t^{2d} k^{2\ell} (2t^d k^\ell p)^{\lceil(t+1)/2\rceil} \leq M^{2d+1} k^{2\ell} (2M^d k^\ell p)^{\lceil(m+2)/2\rceil},$$

so

$$\sum_{t=1}^{M} a(0,t) \leq 2M^{3d+1} k^{3\ell} p,$$

and

$$\sum_{t=m+1}^{M} a(m,t) \leq p \qquad \text{if } m \geq 1$$

and $p$ is sufficiently small. Thus,

$$\sum_{m=1}^{M} (2m+1)^{d+\ell} \sum_{t=m+1}^{M} a(m,t) \leq M(2M+1)^{d+\ell} p.$$

On the other hand, the second bound gives

$$\sum_{m=0}^{R} (2m+1)^{d+\ell} \sum_{t=M+1}^{R} a(m,t) \leq R(2R+1)^{d+\ell} \sum_{t=M+1}^{R} dt \exp(-\delta t)$$

$$\leq R(2R+1)^{d+\ell} dR^2 \exp(-\delta M) \leq p^2.$$

Hence,

$$\mathbb{E}(|\Gamma_G(A,R,x)|) \leq 2M^{3d+1} k^{3\ell} p + M(2M+1)^{d+\ell} p + p^2 \leq 3M^{3d+\ell+1} k^{3\ell} p$$

if $p$ is sufficiently small, as required. $\square$

We need to recall one more easy lemma from [12].

LEMMA 37. *Let $A \subset C([n]^d \times [k]^\ell, r)$. Then for every $1 \leq L \leq \operatorname{diam}([A])$, there exists a connected set $X$ which is internally filled, that is, $X \subset [A \cap X]$, with*

$$L \leq \operatorname{diam}(X) \leq 2L.$$

PROOF. Add newly infected sites one by one, and note that in each step the largest diameter of a component in $[A]$ may jump from at most $L-1$ to at most $2L-1$. Thus, at some point in the process the required set $X$ must appear as a component. $\square$

We are now ready to prove Theorem 1.

PROOF OF THEOREM 1. The upper bound in Theorem 1 was proved in Section 3; the lower bound is an immediate consequence of the following statement. Let $\varepsilon > 0$, and $n \in \mathbb{N}$ be sufficiently large. We shall show that if

$$p = \left( \frac{\lambda(3,3) - \varepsilon}{\log \log n} \right)$$

and $A \in \mathrm{Bin}(B([n]^3, 3), p)$, then

$$\mathbb{P}(\mathrm{diam}([A]) \geq \log n) \leq n^{-30}.$$

Indeed, let $B$ and $k_0$ be given by Corollary 14, let $k \geq k_0$, let $A \in \mathrm{Bin}(B([n]^3, 3), p)$, where $p$ is as above, and suppose $\mathrm{diam}([A]) \geq \log n$. Then, by Lemma 37, there exists an internally filled, connected set $X$ with

$$\frac{\log n - 1}{2} \leq \mathrm{diam}(X) \leq \log n - 1.$$

Let $u, v \in X$ be vertices with $\|u - v\|_\infty + 1 = N := \mathrm{diam}(X)$, and let $R$ be an $[N]^3$ cube, containing $X$, and with $u$ and $v$ on opposite faces of $R$. Write $(x, y, z)$ for an arbitrary element of $R$, where $x, y, z \in [N]$, and let $u \in \{(x, y, z) \in X : x = 1\}$ and $v \in \{(x, y, z) \in X : x = N\}$.

Now, let $m = \lfloor N/k \rfloor$, and partition the cube $R$ into blocks $B_1, \dots, B_m$, each of size $[N]^2 \times [k]$. To be precise, let $B_j = \{(x, y, z) \in R : x \in [(j-1)k + 1, jk]\}$. [If $N$ is not divisible by $k$, then replace $v$ by an element of $\{(x, y, z) \in X : x = km\}$ (the "right-hand face" of $B_m$) and assume that $\{(x, y, z) \in R : x > km\} \subset A$.] Observe that, by our choice of $u$ and $v$, there exists a path in $[A \cap R]$ from the set $\{(x, y, z) \in R : x = 1\}$ to the set $\{(x, y, z) \in R : x = km\}$. We shall use Lemma 35 to show that this is rather unlikely.

Indeed, to do so, we use the following coupling. Replace the thresholds in each block $B_j$ with those of $C([N]^2 \times [k], 2)$, and allow percolation to occur independently in each block. We obtain a set $\bigcup_j [A \cap B_j]$ of infected sites, which we shall denote $\{A\}$. The following claim shows that this is indeed a coupling.

CLAIM 1.   $\{A\} \supset [A \cap R]$.

PROOF. The claim follows easily from the observation that each vertex of $B_j$ has at most one neighbor in $R \setminus B_j$, and internal vertices of $B_j$ [those with $x \notin \{(j-1)k + 1, jk\}$] have no neighbors outside. Indeed, recall that a vertex $w = (x, y, z)$ in $B_j$ originally had threshold 3, and now (in the coupled system) has threshold $2 + I[x \notin \{(j-1)k + 1, jk\}]$. Thus, the threshold of no vertex has increased, and the threshold of those vertices which have a neighbor in $R$ outside $B_j$ have decreased by one. Thus, $\{A\} \supset [A \cap R]$, as claimed. □

For each $j \in [m]$, let $\{A\}(j) = \{A\} \cap B_j$. Now, let $S = [N]^2$, and for each $j \in [m]$, define a two-colored graph $G_j$ on $S \times [2]$ by

$$xy \in E(G_j) \text{ if and only if } \tilde{x} \text{ and } \tilde{y} \text{ are in the same component of } \{A\}(j),$$

where $\tilde{x}$ is the element of $\{(j-1)k+1, jk\} \times [N]^2$ corresponding to $x$ in the natural isomorphism, and

$$x \sim y \Leftrightarrow \text{ there exists an internally filled connected component}$$
$$X \subset \{A\}(j) \text{ such that } x, y \in X \text{ and } \operatorname{diam}(X) \leq B \log n,$$

where $x \sim y$ means $xy$ is a "good" edge, as before, and $B > 0$ was chosen above. Note that $G_j$ is admissible, since $x \sim y$ and $y \sim z$ in $G_j$ implies that $x$ and $z$ are in the same component of $\{A\}(j)$, and so either $x \sim z$, or $xz$ is a bad edge. Note also that the event $x \sim y$ is increasing. We claim that the (random) sequence of admissible two-colored graphs $\mathcal{P} := (G_1, \ldots, G_m) \in \Omega(S, m)$ satisfies the conditions of Lemma 35.

Indeed, recall that $N \leq \log n$, so

$$p \leq \left( \frac{\lambda(3,3) - \varepsilon}{\log N} \right).$$

Choose $0 < \alpha \leq e^{-100k}$, let $\varepsilon' = \varepsilon/3$, and choose $\delta = \delta(\alpha, \varepsilon') > 0$ using Lemma 35. By Corollary 14 (and our choice of $B$ and $k$), for each $j \in [m]$ we have

$$\mathbb{P}(\operatorname{diam}(\{A\}(j)) > B \log N) \leq N^{-\varepsilon},$$

and by Lemma 36, applied with $d = 2$ and $\ell = 1$, for any $v \in V(G_j)$,

$$\mathbb{E}(d_{\mathcal{P}}(v)) = \mathbb{E}(|\Gamma_G(A, B \log N, v)|) \leq \delta$$

if $n$ is chosen to be sufficiently large (and hence $p$ sufficiently small).

Now, conditions (c) and (d) of Lemma 35 are satisfied (for $\delta$ and $\varepsilon'$), by the comments above. Condition (a) is satisfied by construction. Condition (b) follows because if $x \sim y$ and $x' \sim y'$, and there are no bad edges, then either all four points are in the same internally spanned component with diameter at most $B \log n$, or they are in different components of $\{A\}(j)$. So, if $x \not\sim x'$, then the events $x \sim y$ and $x' \sim y'$ must occur disjointly, and so we can apply the van den Berg–Kesten Lemma.

Thus, we may apply Lemma 35, and deduce that

$$\mathbb{P}(X(\mathcal{P})) \leq \alpha^{\lfloor N/k \rfloor} N^2 \leq n^{-40}$$

by our choice of $\alpha$. Summing over all possible rectangles $R$, we see that

$$\mathbb{P}(\operatorname{diam}([A]) \geq \log n) \leq n^{-30},$$

as required. This completes the proof of Theorem 1. $\square$

**Acknowledgments.** Work partly done while R. Morris was visiting the Instituto Nacional de Matemática Pura e Aplicada, Rio de Janeiro, Brazil. The third author would like to thank the Instituto Nacional de Matemática Pura e Aplicada, Rio de Janeiro, Brazil, where a large portion of this research was carried out.

J. Balogh
Department of Mathematics
University of Illinois
1409 W. Green Street
Urbana, Illinois 61801
USA
E-mail: jobal@math.uiuc.edu

B. Bollobás
Trinity College
Cambridge CB2 1TQ
England
and
Department of Mathematical Sciences
University of Memphis
Memphis, Tennessee 38152
USA
E-mail: B.Bollobas@dpmms.cam.ac.uk

R. Morris
Murray Edwards College
University of Cambridge
Cambridge CB3 0DF
England
E-mail: rdm30@cam.ac.uk